\documentclass[12pt,a4paper]{article}
\usepackage{amsmath,amsfonts,amsthm,amsopn,color}
\renewcommand{\ge}{\varepsilon}
\newcommand{\R}{\mathbb{R}}

\newcommand{\RN}{{\mathbb{R}^N}}
\newcommand{\de}{\partial}

\DeclareMathOperator{\cat}{cat}
\DeclareMathOperator{\dv}{div}

\let\eps=\ge
\renewcommand{\a }{\alpha }
\renewcommand{\b }{\beta }
\renewcommand{\d }{\delta }

\renewcommand{\t}{\theta}
\renewcommand{\O}{\Omega}
\newcommand{\G}{\Gamma}
\renewcommand{\S}{\Sigma}
\renewcommand{\L}{\Lambda}
\newcommand{\cl}{{\cal L}}
\newcommand{\Ne}{\mathcal{N}}
\newcommand{\N}{\mathbb{N}}

\newcommand{\tmem}[1]{{\em #1\/}}
\newcommand{\tmop}[1]{\operatorname{#1}}

\newtheorem{theorem}{Theorem}[section]
\newtheorem{lemma}[theorem]{Lemma}
\newtheorem{definition}[theorem]{Definition}
\newtheorem{proposition}[theorem]{Proposition}
\newtheorem{remark}[theorem]{Remark}
\newtheorem{corollary}[theorem]{Corollary}
\renewenvironment{proof}{\noindent{\textbf{Proof\quad}}}{$\hfill\square$\vspace{0.2 cm}\\}
\newenvironment{proof1}{\noindent{\textbf{Proof of Theorem \ref{thm1}\quad}}}{$\hfill\square$\vspace{0.2 cm}\\}
\newenvironment{proof2}{\noindent{\textbf{Proof of Corollary \ref{cor}\quad}}}{$\hfill\square$\vspace{0.2 cm}\\}
\newenvironment{proof3}{\noindent{\textbf{Proof of Theorem \ref{th:multi}\quad}}}{$\hfill\square$\vspace{0.2 cm}\\}
\newenvironment{proofnec}{\noindent{\textbf{Proof of Theorem \ref{thm2}\quad}}}{$\hfill\square$\vspace{0.2 cm}\\}
\newenvironment{proofex}{\noindent{\textbf{Proof of Theorem \ref{th:ex}\quad}}}{$\hfill\square$\vspace{0.2 cm}\\}

\begin{document}

\title{On a class of singularly perturbed elliptic equations in divergence form: existence and multiplicity results}
\author{A.~Pomponio\thanks{Supported by MIUR, national project \textit{Variational methods and nonlinear differential equations} } 
\\ 
SISSA, via Beirut 2/4 \\ I-34014 Trieste 
\\ 
{\it pomponio@sissa.it}
\and S.~Secchi\thanks{Supported by MIUR, national project \textit{Variational methods and nonlinear differential equations} } 
\\ 
Universit\`a di Pisa, via F.~Buonarroti 2 \\ I-56127 Pisa
\\ 
{\it secchi@mail.dm.unipi.it}}

\date{ }

\maketitle

\section{Introduction}

The aim of this paper is to study the existence and the concentrating
behavior of solutions to the following problem:

\begin{equation}\label{eq:E1}
\begin{cases}
-\ge^2\dv \big( J(x)\nabla u \big) + V(x)u=u^p & \textrm{ in } \RN
\\
u>0 & \textrm{ in } \RN
\\
\lim_{|x|\to \infty} u(x)=0
\end{cases}
\end{equation}
where $N \geq 3$, $p\in \left(1, \frac{N+2}{N-2} \right)$, $V\colon \RN \to \R$, 
$J\colon\RN \to \R^{N \times N}$ are $C^1$ functions. Here the symbol $\R^{N\times N}$
stands for the set of ($N\times N$) real matrices.

Such a problem, at least in the case $J\equiv I$, where $I$ is the identity matrix in $\R^{N \times N}$, 
arises naturally when seeking \textit{standing waves} of the nonlinear Schr\"odinger 
equation with bounded potential $V$, that are solutions of the form
\[
\psi(t,x)={\rm e}^{i\hbar^{-1}t}u(x)
\]
of the following problem
\[
i\hbar \frac{\partial\psi}{\partial t}=-\hbar^2\varDelta\psi +V(x)\psi-|\psi|^{p-1}\psi, \quad x\in\RN,
\]
where $\hbar$ denotes the Planck constant, $i$ is the unit imaginary. The usual strategy is to put $\ge=\hbar$ and then to study what happens when $\ge \to 0$.

Problem \eqref{eq:E1}, at least with $J\equiv I$, has been studied extensively in several works, 
see e.g. \cite{ABC,dPF,dPF2,dPF4,FW,Li,O,Oh,Rab,W,WZ}. 

It is known that if \eqref{eq:E1} has a solution concentrating at some $z_0$, then $\nabla V(z_0)=0$. 
Conversely, if $z_0$ is a critical point of $V$ with some stability properties, then \eqref{eq:E1} 
has a solution concentrating at $z_0$ (see for example \cite{ABC,AMS,Li}).

Two main strategies have been followed. A first one, initiated by Floer and Weinstein \cite{FW}, relies on 
a finite dimensional reduction. The second one has been introduced by del Pino and Felmer \cite{dPF} 
and is based on a penalization technique jointly with local linking theorems. 

In the present paper, we study \eqref{eq:E1} in the case $J\not\equiv I$. Our research is motivated 
by \cite{S}, where a general class of singularly perturbed quasilinear equation on $\RN$,
\begin{equation}\label{eq:Squa}
-\eps^2 \dv \big( J(x,u)\nabla u \big) + 
\frac{\eps^2}{2}\langle D_s J(x,u) \nabla u \mid \nabla u\rangle
+V(x)u = f(u),
\end{equation}  
is studied by means of non-smooth critical points theory. If $J$ depends only on $x$ and $f(u)=u^p$, 
then \eqref{eq:Squa} becomes \eqref{eq:E1}.

We observe that it is in general impossible to reduce the second-order operator in equation \eqref{eq:E1} to the standard Laplace operator in the whole $\RN$ by means of a single change of coordinates. This phenomenon especially appears in high dimension $N > 3$, as already remarked in \cite{CH}, chapter III.

On $V$ and $J$ we will make the following assumptions:
\begin{description}
\item[(V)] $V\in C^1 (\RN, \R)$ and $\inf_{\RN} V = \alpha >0$;
\item[(J)] $J\in C^1 (\RN, \R^{N\times N})$, $J$ is bounded; moreover, $J(x)$ is, for each $x\in\RN$, 
a symmetric matrix, and
\begin{equation}\label{eq:ell}
(\exists \nu >0 )(\forall x\in\RN)(\forall\xi\in\RN\setminus \{0\}):\;
\left\langle J(x)\xi,\xi\right\rangle \geq \nu |\xi|^2.
\end{equation}
\end{description}

Let us introduce an auxiliary function which will play a crucial r\^ole in the study of \eqref{eq:E1}. 
Let $\G\colon \RN \to \R$ be a function so defined:
\begin{equation}\label{eq:Gamma}
\G(z)=V(z)^{\frac{p+1}{p-1}- \frac{N}{2}}
\left(\det J(z)  \right)^{\frac{1}{2}}.
\end{equation}
Let us observe that by {\bf (J)}, $\G$ is well defined.

We now state the main results of this work. We will see that $\G$ gives, 
roughly speaking, a sufficient condition and a necessary one to have concentrating solutions around a point.

\begin{theorem}\label{th:1}
Suppose {\bf (V)} and {\bf (J)} hold. 
Suppose that there exists a compact domain $\L \subset \RN$ such that 
\begin{equation*}
\min_\L \G<\min_{\de \L}\G.
\end{equation*}
Then, for all $\eps >0$ sufficiently small, there exists a solution
$u_\eps \in H(\RN) \cap C(\RN)$ of \eqref{eq:E1} with $\int_\RN V(x) u^2 d x<+\infty$. 
This solution has only one global maximum point $x_\eps\in \RN$ and we have that 
$\G (x_\ge)\to \min_\Lambda \G$ as $\eps \to 0$ and 
\[
\lim_{\ge\to 0} u_\ge (x) =0 \mbox{\quad for all $x\neq x_\ge$}.
\]
\end{theorem}

\begin{theorem}\label{th:2}
Assume, in addition to assumptions {\bf (V)} and {\bf (J)}, 
that $V$ is bounded from above. Let $\{u_{\ge_j}\}$ be a
sequence of solutions of \eqref{eq:E1} such that for all $\ge >0$ there
exist $\rho >0$ and $j_0>0$ such that for all $j \geq j_0$ and for all points
$x$ with $|x-z_0| \geq \ge_j \rho$, there results
\[
u_{\ge_j}(x) \leq \ge.
\]
Then $z_0$ is a critical point of $\G$.
\end{theorem}

Theorem \ref{th:1} and \ref{th:2} are simplified version of, respectively, 
Theorem \ref{thm1} and \ref{thm2}, which hold under weaker assumptions. 
Indeed we can treat \eqref{eq:E1} in a more general domain, instead of $\RN$, and 
with a more general nonlinearity.
More precisely, let us consider:
\begin{equation}\label{eq:E}
\begin{cases}
-\ge^2\dv \big( J(x)\nabla u \big) + V(x)u=f(u) & \textrm{ in } \Omega
\\
u>0 & \textrm{ in } \O
\\
u=0 & \textrm{ on }  \de\O
\end{cases}
\end{equation}
where $\O$ is an open domain of $\RN$, possibly unbounded, and $f \colon \R^+ \to \R$ is a $C^1$ such that:
\begin{description}
\item[(f1)]$f(u)=o(u)$ as $u \to 0^+$;
\item[(f2)]for some $p \in \left(1, \frac{N+2}{N-2} \right)$ there holds 
\[
\lim_{u \to +\infty} \frac{f(u)}{u^p}=0 ;
\]
\item[(f3)]for some $\t \in (2, p+1)$ we have
\[
0< \t F(u)\leq f(u)u \quad \textrm{ for all }u >0,
\]
where $F(u)= \int_0^u f(t) d t$;
\item[(f4)]the function 
\[
u \in (0,+\infty) \mapsto \frac{f(u)}{u}
\]
is increasing.
\end{description}

The proof of Theorem \ref{th:1} is based on the penalization technique used in 
\cite{dPF}, adapted to our case. See Section 2 and 3.

If $z_0$ is a common minimum point of $V$ and $J$, $J$ depends only by $x$ and $f(u)=u^p$, 
then our result becomes a particular case of Theorem 1.1 of \cite{S}. On the other side, 
\cite{S} considers the case when $V$ and $J$ have a common minimum point, only.

In Section 4 we prove Theorem \ref{th:2} using a recent version of Pucci--Serrin variational 
identity (see \cite{DMS}).

In Section 5 we consider \eqref{eq:E1} assuming that $V$ and $J$ satisfy, in addition to 
hypotheses {\bf (V)} and {\bf (J)}: 
\begin{description}
\item[(V1)] $V\in C^{2}(\RN, \R)$, $V$ and $D^{2}V$ are bounded;
\item[(J1)] $J\in C^{2}(\RN,\R^{N\times N})$, $J$ and $D^{2}J$ are bounded.
\end{description}
Then the following theorem holds.
\begin{theorem}\label{th:ex}
Let {\bf (V-V1)} and {\bf (J-J1)} hold. Then for $\eps>0$ small, (\ref{eq:E1}) has a 
solution concentrating in $z_0$, provided that one of the two following conditions holds:
\begin{description}
\item[$(a)$] $z_0$ an isolated local strict minimum or maximum of $\G$;
\item[$(b)$] $z_0$ is a non-degenerate critical point of $\G$.
\end{description}
\end{theorem}

The proof of Theorem \eqref{th:ex} relies on a finite dimensional reduction, precisely on 
the perturbation technique developed in \cite{AB,ABC,AMS}.

The last section is devoted to the proof of some multiplicity results, see Theorem \ref{th:multi}, \ref{th:main} 
and \ref{th:CL}. As before we distinguish between two cases, according to the use of penalization 
method of \cite{dPF} or of the perturbation method of \cite{AMS}.

All the results contained in Section 5 and 6 seem to be new and do not appear in \cite{S}.

\bigskip

\noindent {\bf Acknowledgments}\quad The authors wish to thank Professor Antonio Ambrosetti
 for suggesting the problem and for useful discussions. S.~S. thanks also SISSA for the kind hospitality.

\bigskip

\begin{center}
\textbf{Notation}
\end{center}

\begin{itemize}
\item If $F$ is a $C^1$ map on a Hilbert space $H$, we denote by $DF(u)$ its Fr\'{e}chet derivative at $u\in H$.
\item For $x,y\in\RN$, we denote by $\langle x\mid y \rangle$ the ordinary inner product of $x$ and $y$.
\item $C$ denotes a generic positive constant, which may also vary from line to line.
\item $o_h (1)$ denotes a function that tends to $0$ as $h\to 0$.
\end{itemize}

\section{The ground-energy function $\S$}\label{sec:S}

In this section we present a more general version of Theorem \ref{th:1} and we  
introduce the ground-energy function $\S$ that has a crucial r\^ole in the sequel and that, 
at least when $\O=\RN$ and $f(u)=u^p$, is equal to $\G$ up to a  constant factor. 

We work in the weighted space 
\[
H_V(\O)=\left\{ u\in H^1_0(\O) \colon \int_\O V(x) u^2 d x <+\infty 
\right\},
\]
endowed with the norm 
\[
\|u\|\sp 2 = \int_\O \left( |\nabla u|^2 + V(x)u^2 \right) dx.
\]

Our assumptions imply that the functional $I^\eps \colon H_V(\O) \to \R$
defined by
\begin{equation}\label{eq:I-eps}
I^\eps (u)= \frac{\eps^2}{2}\int_\O \left\langle J(x)\nabla u\mid \nabla
u\right\rangle + 
\frac{1}{2}\int_\O V(x) |u|^2 - \int_\O F(u)
\end{equation}
is of class $C^2(\O)$. Moreover, \eqref{eq:E} is the Euler-Lagrange
equation associated to $I^\eps$, so that we will find solutions of
\eqref{eq:E} as critical points of $I^\eps$.

We define the {\it ground-energy function} $\S(z)$ as the ground energy associated with
\begin{equation}\label{eq:Ez}
-\dv \big( J(z)\nabla u \big) + V(z)u=f(u) \quad \textrm{ in } \RN,
\end{equation}
where $z\in \RN$ is seen as a (fixed) parameter. More precisely,
\eqref{eq:Ez} is associated to the functional defined on $H^1(\RN)$
\begin{equation}\label{eq:E-z}
I_z (u)= \frac{1}{2}\int_\RN \left\langle J(z)\nabla u\mid \nabla
u\right\rangle d x + 
\frac{1}{2}\int_\RN V(z) |u|^2 d x- \int_\RN F(u) d x.
\end{equation}
If $\mathcal{N}_z$ is the \textit{Nehari manifold} of
\eqref{eq:E-z}, that is
\[
\mathcal{N}_z = \left\{u \in H^1(\RN) \mid \mbox{$u\neq 0$ and $D I_z (u)[u]=0$}
\right\},
\]
we have by definition
\begin{equation}\label{eq:Xi}
\S(z)=\inf_{u \in \mathcal{N}_z}I_z (u).
\end{equation}

\begin{remark}\label{rem:SG}
As already said, when $\O=\RN$ and $f(u)=u^p$, we have that there exists a positive 
constant $C>0$ such that
\[
\S(z)=C \G(z).
\]
Indeed, if $U$ is the unique radial solution in $H^1(\RN)$ of
\begin{equation*}
\begin{cases}
-\varDelta U+ U = U^p \quad \textrm{ in } \RN 
\\
U>0,
\end{cases}
\end{equation*} 
then it is easy to see that
\[
\S(z)= \left( \frac{1}{2}- \frac{1}{p+1}\right)
V(z)^{\frac{p+1}{p-1}- \frac{N}{2}}
\left(\det J(z)  \right)^{\frac{1}{2}}
\int_\RN U^{p+1}.
\]
\end{remark}

It is easy to see that $\Ne_z \neq \emptyset$ and moreover the following lemma holds.

\begin{lemma}\label{le:retta}
For all $u \in H^1(\RN)$ such that $u$ is positive on a set of positive measure, 
there exists a unique maximum $t(u)>0$ of
\[
\phi \colon t\in (0, +\infty) \mapsto I_z(t u).
\] 
In particular, $t(u)u \in \mathcal{N}_z$.
\end{lemma}

\begin{proof}
Let us observe that if $\phi'(t)=0$, then
\[
\int_\RN \langle J(z) \nabla u \mid \nabla u \rangle + V(z) u^2 = \int_\RN \frac{f(t u)u}{t}\, dx
\]
and so, by \textbf{(f4)}, $\phi$ has at most one critical value. By \textbf{(f1-2)}, 
$I_z(0)=0$, $D I_z(0)=0$ and $D^2I_z(0)$ is strictly positive-definite in a neighborhood of $0$ and so $\phi(t)>0$ for 
$t$ small. Moreover, since
\[
\int_\RN F(t u)=\int_{\{x\in\RN \colon u(x)>0\}}F(t u),
\]
by \textbf{(f3)} there results $\phi(t)<0$ for big $t$'s.
\end{proof}

The following proposition gives us some useful properties of $\S$ (see also \cite{WZ}).

\begin{proposition}\label{prop:1}
Let the assumptions {\bf (V)}, {\bf (J)}, {\bf (f1-4)} hold. Then:
\begin{enumerate}
    \item the map $\S$ is well-defined and locally lipschitz;
    \item the partial derivatives, from the left and the right, of $\S$
exist at every point, and moreover
\begin{equation*}
    \left( \frac{\de\S}{\de s_i}\right) \sp l = \sup_{v_s\in S^s}
\left[ \frac{1}{2}\int_\RN \langle \frac{\de J}{\de s_i} \nabla v_s \mid \nabla v_s
\rangle + 
\frac{1}{2} \frac{\de V}{\de s_i} \int_\RN |v_s|^2
\right]
\end{equation*}
\begin{equation*}
    \left( \frac{\de\S}{\de s_i}\right) \sp r = \inf_{v_s\in S^s}
\left[ \frac{1}{2} \int_\RN \langle \frac{\de J}{\de s_i} \nabla v_s \mid \nabla v_s
\rangle +
\frac{1}{2} \frac{\de V}{\de s_i} \int_\RN |v_s|^2
\right]
\end{equation*}
where $S^s$ is the set of ground states corresponding to the energy
level $\S (s)$.
\end{enumerate}    
\end{proposition}

\begin{proof}
First of all, the set $S^s$ is non-empty. Indeed, since $s$ is fixed, 
we can find a matrix $T=T(s)\in \operatorname{GL}(n)$ such that
\[
T^t J(s) T =I \mbox{\quad (the identity matrix of order $n$)}.
\]
By the change of variables $x \mapsto Tx$, the equation
\[
-\dv \left(J(s)\nabla v \right) + V(s) v = f(v) \quad \textrm{ in } \RN
\]
can be reduced to
\begin{equation}
\label{eq:rescaled}
-\varDelta U + V(s)U=f(U) \quad \textrm{ in } \RN.
\end{equation}
This change of variables rescales the functional $I_z$ by a constant. Since it is well 
known that equation \eqref{eq:rescaled} has a ground state solution, for
each $s\in\RN$, it immediately follows that $S^s \neq \emptyset$.

Let us observe that if $v_t \in \Ne_t$, since it satisfies
\[
\int_\RN \langle J(t) \nabla v_t \mid \nabla v_t \rangle + V(t) |v_t|^2 = \int_\RN f(v_t)v_t \, dx,
\]
$v_t>0$ on a set of positive measure and so we can apply the Lemma \ref{le:retta}.
Therefore, given $s$, $t\in \RN$,
there exists precisely one positive number $\theta (s,t)$ such that $\theta (s,t)v_t\in\mathcal{N}_s$. 
By definition, this means that
\begin{equation*}
    \int_\RN \langle J(s)\nabla v_t \mid \nabla v_t \rangle + \int_\RN
V(s)|v_t|^2 = \int_\RN \frac{f(\theta (s,t)v_t)v_t}{\theta (s,t)}.
\end{equation*}
Moreover, $\theta (t,t)=1$. Collecting these facts, we see that, by the
implicit function theorem, $\t$ is differentiable with respect to the first
variable. From its very definition, $\theta (s,t)$ is bounded for $s$
and $t$ bounded in $\RN$. Let us now observe that
\begin{multline*}
    I_s(\theta (s,t)v_t) = \frac{\theta (s,t)^2}{2}  \int_\RN \langle
J(s)\nabla v_t \mid \nabla v_t \rangle +
\\ 
+ \frac{\theta(s,t)^2}{2} \int_\RN V(s)|v_t|^2 
    - \int_\RN F(\theta (s,t)v_t).
\end{multline*} 
The gradient of the function $s\mapsto I_s(\theta (s,t)v_t)$ is thus
\begin{multline*}
    \frac{\de}{\de s}I_s(\theta (s,t)v_t)=\frac{\theta(s,t)^2}{2} \int_\RN \langle
\nabla J(s)\nabla v_t\mid \nabla v_t \rangle +
\frac{\theta(s,t)^2}{2} \nabla V(s) \int_\RN |v_t|^2
\\
    +\theta (s,t) \frac{\de \theta}{\de s} \left( \int_\RN \langle
J(s)\nabla v_t\mid \nabla v_t  \rangle+ V(s) \int_\RN |v_t|^2 \right) -
\int_\RN f(\theta (s,t)v_t) v_t \frac{\de \theta}{\de s} \\
= \frac{\theta (s,t)^2}{2} \int_\RN \langle \nabla J(s) \nabla v_t \mid \nabla v_t
\rangle +
\frac{\theta (s,t)^2}{2} \nabla V(s) \int_\RN |v_t|^2
\end{multline*}
because $\theta (s,t)v_t \in \mathcal{N} \sb s$. From this representation,
the mean value theorem and the local boundedness of $\theta$,
it follows that for all $R>0$ there exists $M>0$ such that for all
$s_1$ and $s_2$ with $|s_1|<R$, $|s_2|<R$:
\[
|\S (s_1) - \S (s_2) | \leq M |s_1-s_2|.
\]
This proves the first statement. The proof of the second statement can
be modeled on the similar results of \cite{WZ} together with
\cite{RS}. We omit the details.
\end{proof}

\begin{remark}
Some uniqueness conditions for the limiting equation appear in \cite{ST}.
\end{remark}

The next step is the proof of a fundamental equality between the ground state level and 
the mountain-pass one. This kind of result is well known at least in the case $J$ equal 
to the identity matrix (see, for example, \cite{Rab}).

\begin{proposition}\label{prop:ugua}
\begin{equation*}
\inf_{\gamma\in\mathcal{P}_z} \max_{0\leq t\leq 1} I_z(\gamma(t)) =
\inf_{u\in \mathcal{N}_z} I_z(u) = \S (z),
\end{equation*}
where
\[
\mathcal{P}_z = \left\{ \gamma\in C\left([0,1],H^1(\RN)\right) \colon \gamma
(0)=0 \;{\rm and}\; I_z(\gamma(1))<0\right\}.
\]
\end{proposition}

\begin{proof}
We now show that for all path $\gamma\in\mathcal{P}_z$, there exists $t_0>0$ 
such that $\gamma (t_0) \in \mathcal{N}_z$.

First of all, let us observe that by our assumptions of $f$, together with the ellipticity
of $J$ and the definition of $\alpha = \inf V >0$ we get
\begin{eqnarray*}
DI_z (u)[u] &=& \int_\RN \langle J(z) \nabla u \mid \nabla u \rangle +
V(z) u^2 - \int_\RN f(u)u\, dx \\
&\geq& \nu \int_\RN |\nabla u|^2 + \alpha \int_\RN |u|^2 -c\int_\RN
|u|^{\t+1} \\
&\geq& \min \{ \nu,\alpha \} \|u\|_{H_0^1(\Omega)}^2 - c
\|u\|_{L^{\t+1}(\Omega)}.
\end{eqnarray*}
Therefore, fix a path $\gamma\in\mathcal{P}_z$ joining $0$ to some $v\neq 0$ such that
$I_z (v)<0$. Hence
\[
DI_z(\gamma (t))[\gamma(t)] > 0
\]
for $t>0$ small enough. On the other hand, since $v\neq 0$, we have
\[
DI_z(v)[v] < 2 I_z(v) \leq 0.
\]
By the intermediate value theorem, there exists $t_0\in [0,1]$ such that
\[
\gamma (t_0) \in \mathcal{N}_z.
\]
This shows also that
\begin{equation}
\label{eq:minimax}
\inf_{\gamma\in\mathcal{P}_z} \max_{0\leq t\leq 1} I_z(\gamma(t)) \geq
\inf_{u\in \mathcal{N}_z} I_z(u) = \S (z).
\end{equation}
But equality actually holds in \eqref{eq:minimax}. Indeed, by Lemma \ref{le:retta},
\begin{equation}\label{eq:max}
u\in\mathcal{N}_z \quad \textrm{ if and only if } \quad I_z (u) = \max_{\tau \geq 0} I_z (\tau u).
\end{equation}
It then follows that the mountain-pass level corresponds to the least
energy among the energies of all solutions, that is
\begin{equation*}
\inf_{\gamma\in\mathcal{P}_z} \max_{0\leq t\leq 1} I_z(\gamma(t)) =
\inf_{u\in \mathcal{N}_z} I_z(u) = \S (z).
\end{equation*}
\end{proof}

\begin{remark}
It is shown in \cite{JT} that the level of mountain--pass of $I_z$ coincides with the infimum of the energies among all critical points of $I_z$ under less stringent assumptions on $f$. Anyway, to prove the regularity of $\Sigma$ we are not able to weaken our set of hypotheses {\bf (f1--4)}.
\end{remark}

Our main result about existence for \eqref{eq:E} is the following theorem.

\begin{theorem}\label{thm1}
Suppose {\bf (V)}, {\bf (J)} and {\bf (f1-4)} hold. Suppose that there
exists a compact domain $\L \subset \O$ such that 
\begin{equation}\label{eq:min}
\min_\L \S<\min_{\de \L}\S.
\end{equation}
Then, for all $\eps >0$ sufficiently small, there exists a solution
$u_\eps \in H_V(\O) \cap C(\O)$ 
of \eqref{eq:E}. Moreover, 
this solution has only one global maximum point $x_\ge\in\Lambda$ and we have that 
$\Sigma (x_\ge)\to \min_\Lambda \Sigma$ as $\eps \to 0$  and
\[
\lim_{\ge\to 0} u_\ge (x) =0 \mbox{\quad for all $x\neq x_\ge$}.
\]
\end{theorem}

\begin{remark}
As noticed in \cite{S}, we have that
\[
\lim_{\eps \to 0}\|u_\eps \|=0.
\]
\end{remark}

The next corollary shows that results of \cite{S}, at least in the case of our differential operator, 
are a particular case of Theorem \ref{thm1}. More precisely, as said in the introduction, we will prove 
that if $J$ and $V$ have a local strict minimum in $z_0$, then $z_0$ is also a local strict minimum 
for $\S$ and so we can apply Theorem \ref{thm1}.

\begin{corollary}\label{cor}
Suppose {\bf (V)}, {\bf (J)} and {\bf (f1-4)} hold. Suppose that there
exist a compact domain $\L \subset \O$ and a $z_0 \in \L$ which is a minimum point 
for $V$ and $J$ in $\L$ and a strict minimum point for $V$ (resp. $J$), 
in the sense that 
\begin{equation}\label{1}
V(z_0)<\min_{\de \L}V \quad \left({\textrm resp. }\,   V(z_0)\leq\min_{\de \L}V \right) 
\end{equation}
and, for all $\xi \in \RN \setminus \{ 0 \}$ and for all $z \in \de \L$, 
\begin{equation}\label{2}
\langle J(z_0) \xi \mid \xi \rangle \leq \langle J(z) \xi \mid \xi \rangle \quad
\left({\textrm resp. }\,   \langle J(z_0) \xi \mid \xi \rangle < \langle J(z) \xi \mid \xi \rangle \right).
\end{equation}
Then \eqref{eq:min} holds and hence the conclusion of Theorem \ref{thm1} continues to be true.
\end{corollary}

We will prove Theorem \ref{thm1} and Corollary \ref{cor} in the next section.

\section{The penalization scheme}\label{sec:pen}

Since the domain $\Omega$ is in general unbounded, a direct application of
critical point theory does not, as a rule, provide a solution to
\eqref{eq:E}. Although the functional $I^\eps$ has a good geometric
structure, it does not satisfy the Palais-Smale condition. Thus, as a
first step, we replace $I^\eps$ with a different functional that
satisfies (PS)$_c$ at all levels $c\in\R$, and finally prove that, as
$\ge$ gets smaller, the critical points of this new
functional are actually solutions of \eqref{eq:E}. This technique was
introduced by del Pino and Felmer in \cite{dPF}, and then used by several authors. 
The main advantage of this scheme is
that, unlike the direct application of some {\it Concentration-Compactness argument} 
as in \cite{WZ}, we do not have
to impose any comparison assumption between the values of $\S$ at zero
(say) and at infinity.

Following the scheme of \cite{dPF} (see also \cite{S}), we will define a penalization of the functional
$I^\ge$, which satisfies 
the Palais-Smale condition. Let $\t$ be the number given in {\bf (f3)}.
Let $\ell>0$ be the unique value such that $f(\ell)/\ell=\alpha/k$,
where $\alpha$ is defined in {\bf (V)} and $k>\t/(\t -2)$.

We penalize the nonlinearity $f$ in the following way. Define $\tilde{f}
\colon \R \to \R$ by
\[
\tilde{f} (u)=
\begin{cases}
(\alpha/k) u & \text{if $u > \ell$} \\
f(u) & \text{if $0\leq u \leq \ell$} \\
0 & \text{if $u <0$}.
\end{cases}
\]

We now define $g \colon \O \times \R \to \R$ as
\[
g(x,u)=
\begin{cases}
\chi_\L (x) f(u)+ (1-\chi_\L (x) )\tilde{f} (u) & \textrm{ if } u
\geq 0 \\
0 & \textrm{ if } u <0
\end{cases}
\]
where $\chi_\L$ is the characteristic function of the set $\L$, and let
$G$ be the primitive of $g$, that is
\[
G(x, u)=\int_0^u g(x, \tau)d \tau.
\]
By straightforward calculations, assumptions {\bf (f1-4)} imply:
\begin{description}
\item[(g1)]$g(x,u)=o(u)$ as $u \to 0^+$, uniformly in $x\in \O$;
\item[(g2)]$\lim_{u \to +\infty} g(x,u)/u^{p}=0$ for some 
$p \in \left(1, \frac{N+2}{N-2} \right)$;
\item[(g3-i)]for some $\t \in (2, p+1)$ we have
\[
0 < \t G(x,u)\leq g(x,u)u \quad \textrm{ for all } x\in \L,\, u >0;
\]
\item[(g3-ii)] for some $k>\frac{\t}{\t -2}$ there holds
\[
0\leq 2 G(x,u)\leq g(x,u)u \leq \frac{1}{k}V(x)u^2  \quad \textrm{
for all } x\notin \L,\, u >0;
\]
\item[(g4)]the function $u \mapsto \frac{g(x, u)}{u}$ is
increasing for all $x\in \L$.
\end{description}

The penalized functional will be $E^\eps \colon H_V(\O) \to \R$, where
\begin{equation}\label{eq:E-eps}
E^\eps (u)= \frac{\eps^2}{2}\int_\O \left\langle J(x)\nabla u\mid \nabla
u\right\rangle + \frac{1}{2}\int_\O V(x) |u|^2 - \int_\O G(x,u).
\end{equation}

Under our assumptions $E^\eps$ satisfies the (PS) condition, as we prove
in the next lemma.

\begin{lemma}
Let $\{u_h \}$ be a sequence in $H_V(\O)$ such that $E^\eps (u_h) \to c
\in \R$ 
and $DE^\eps (u_h) \to 0$. Then $\{u_h \}$ has a strongly convergent
subsequence.
\end{lemma}

\begin{proof}
As first step, we show that $\{u_h \}$ is bounded. Since $E^\eps (u_h)
\to c$ 
we have
\begin{multline*}
\frac{\t}{2}\eps^2 \int_\O \left\langle J(x)\nabla u_h\mid \nabla
u_h\right\rangle + 
\frac{\t}{2}\int_\O V(x) |u_h|^2 
\\
\leq \int_\L g(x,u_h)u_h 
+ \frac{\t}{2 k}\int_{\O \setminus \L} V(x) |u_h|^2 +\t c + o(1).
\end{multline*}
Moreover, since $DE^\eps (u_h)[u_h]=o(\|u_h\|)$, 
\[
\eps^2 \int_\O \left\langle J(x)\nabla u_h\mid \nabla u_h\right\rangle + 
\int_\O V(x) |u_h|^2 
\geq \int_\L g(x,u_h)u_h +o(\|u_h\|).
\]
Therefore
\begin{multline*}
\min\left\{ \left(\frac{\t}{2} -1 \right)\eps^2,\,
\frac{\t}{2}-\frac{\t}{2k}-1 \right\}
\int_\O \left\langle J(x)\nabla u_h\mid \nabla u_h\right\rangle + V(x)
|u_h|^2 
\\
\leq \t c + o(1) + o(\|u_h\|) 
\end{multline*}
and so the boundedness of  $\{u_h \}$ follows from \eqref{eq:ell}. 

Up to subsequence, we have that $u_h \to u$ weakly and point-wise almost
everywhere in $\O$. To show that this convergence is actually strong,
it suffices to prove that, for all $\delta >0$, there exists $R>0$
such that 
\[
\limsup_{k \to \infty} \int_{\O \setminus B_R} |\nabla u_h|^2 + V(x)
|u_h|^2 < \delta.
\]
We take $R>0$ so large that $\L \subset B_{R/2}$. Let $\eta_R \in C^2(\O)$ be a function such that, 
$\eta_R =0$ in $B_{R/2}$, $\eta_R =1$ in $\O \setminus B_{R}$, $0\leq
\eta_R \leq 1$ in $\O$ and 
$|\nabla\eta_R| \leq C/R$ in $\O$.
Since $\{u_h \}$ is bounded,
\[
\lim_{k \to \infty} DE^\eps (u_h)[\eta_R u_h]=0.
\]
Therefore
\begin{gather*}
\int_\O \left (\left\langle J(x)\nabla u_h\mid \nabla u_h\right\rangle + 
V(x) |u_h|^2 \right)\eta_R
+ \int_\O \left\langle J(x)\nabla u_h\mid \nabla \eta_R \right\rangle
u_h =
\\
\int_\O g(x,u_h)u_h \eta_R + o(1)
\leq \frac{1}{k} \int_{\O} V(x) u_h^2 \eta_R + o(1),
\end{gather*}
and so
\[
\int_{\O \setminus B_R} |\nabla u_h|^2 + V(x) u_h^2 \leq
\frac{C}{R}\|u_h\|_{L^2}\|\nabla u_h\|_{L^2}
+ o(1).
\]
We conclude the proof by letting $R\to+\infty$.
\end{proof}

Since by Proposition \ref{prop:1} we know that $\S$ is a continuous function, we can assume without
loss of generality that there exists $z_0\in \Lambda$ such that
\[
\S (z_0) = \min_\L \S.
\]
Hence the main assumption of the Theorem \ref{thm1} can be stated as
\[
\S (z_0) < \min_{\de\L} \S.
\]
To save notation, we will write $I_0$ instead of $I_{z_0}$.

Let us set
\begin{equation*}
\bar c=\inf_{\gamma\in \mathcal{P}_0} \max_{t\in [0,1]} I_0 (\gamma(t))=\S(z_0),
\end{equation*}
where
\[
\mathcal{P}_0 = \left\{ \gamma\in C\left([0,1],H^1(\RN)\right) \colon
\hbox{$\gamma (0)=0$ and $I_0 (\gamma (1))<0$}\right\}.
\]

\begin{lemma}\label{le:dis}
For all $\ge$ sufficiently small, there exists a critical point $u_\ge
\in H_V(\Omega)$ of $E^\ge$ such that
\[
E^\ge (u_\ge) \leq \ge^N \bar c + o(\ge^N).
\]
\end{lemma}

\begin{proof}
We already know that $E^\ge$ satisfies the (PS) condition at any level.
By a standard minimax argument over the set of paths
\[
\mathcal{P}_\ge = \left\{ \gamma \in C\left([0,1],H_V(\O)\right) \colon
\hbox{$\gamma (0)=0$ and $E^\ge (\gamma (1))<0$} \right\},
\]
we can find a critical point $u_\ge$ such that
\[
E^\ge (u_\ge)=\inf \sb {\gamma \in \mathcal{P}_\ge} \max \sb{t\in [0,1]}
E^\ge (\gamma (t)).
\]
Since $\bar c$ is a {\it mountain-pass level} of $I_0$, for all $\delta
>0$ there exists a path $\gamma \colon [0,1] \to H^1(\RN)$ such that
\[
\bar c \leq \max_{0\leq t\leq 1} I_0 (\gamma (t)) \leq \bar c + \delta,\quad
\gamma (0)=0,\quad
I_0 (\gamma (1))<0.
\]
Let $\zeta \in C^2(\RN)$ be a cut-off function such that $\zeta =1$ in a
neighborhood of $z_0$. Define a path in $H_V(\O)$ by
\[
\Gamma \sb \ge (\tau) \colon x \mapsto \zeta (x) \gamma (\tau) \left(
\frac{x-z_0}{\ge} \right).
\]
By direct computation,
\begin{multline*}
    E^\ge (\Gamma_\ge (\tau)) = \ge^N \left\{ \frac{1}{2} \int_\RN
\langle J(z_0)\nabla \gamma (\tau)\mid \nabla \gamma (\tau)\rangle +
\right.
\\ 
+ \frac{1}{2} \int_\RN V(z_0)|\gamma (\tau)|^2 
    \left. - \int_\RN F(\gamma(\tau)(\cdot))\right\} + o(\ge^N),
\end{multline*}
that is
\[
E^\ge (\Gamma_\ge (\tau)) = \ge^N I_0 (\gamma(\tau)) + o(\ge^N)
\]
as $\ge\to 0$. But $\Gamma_\ge\in\mathcal{P}_\ge$, so  that
\begin{eqnarray*}
  E^\ge (u_\ge) &=& \inf_{\gamma\in \mathcal{P}_\ge} \max_{t\in [0,1]}
E^\ge (\gamma (t)) \\
   &\leq& \max_{t\in [0,1]} E^\ge (\Gamma_\ge (t)) = \ge^N \max_{t\in
[0,1]} I_0 (\gamma (t)) + o(\ge^N)\\
  &\leq& \ge^N \bar c  + \delta \ge^N + o(\ge^N).
\end{eqnarray*}  
Since $\delta>0$ was arbitrary, the proof is complete.
\end{proof}

\begin{remark}
By the uniform ellipticity of $J$ and standard regularity theorems (see
\cite{G}), the element $u_\ge$ actually belongs to $C(\overline{\L})$.
\end{remark}

The next proposition is somehow the key ingredient.

\begin{proposition}\label{prop:key}
Let $u_\ge\in H_V(\O)$ be the critical point of $E^\eps$ found in the
previous lemma. Then
\[
\lim_{\ge\to 0} \max_{\de\L} u_\ge =0.
\]
Moreover, $u_\eps$ has only one global maximum point $x_\ge\in\Lambda$ and we have that 
$\Sigma (x_\ge)\to \min_\Lambda \Sigma$ as $\eps \to 0$  and
\[
\lim_{\ge\to 0} u_\ge (x) =0 \mbox{\quad for all $x\neq x_\ge$}.
\]
\end{proposition}

\begin{proof}
The proof of this Proposition will be performed in several steps. 
First of all, the following claim implies the Proposition 
except for the uniqueness of the global maximum.
\vspace{0.5 cm}

\noindent{\tmem{Claim 1: If $\{ \varepsilon_h \}$ is a sequence of positive real
numbers
converging to zero, and $\{ x_h \} \subset \Lambda$}} {\tmem{is a
sequence of
points in $\Lambda$ such that}}
\[ u_{\varepsilon_h} ( x_h ) \geq c > 0,  \]
{\tmem{then
\[ \lim_{h \rightarrow \infty } \S ( x_h ) = \min_{\Lambda} \S . \]}}

Indeed, suppose that the statement of the Proposition is false. Then,
up to a subsequence, we may suppose that there exists a sequence $\{x_n \}\subset \de\L$ 
such that $x_h \to \bar{x} \in \partial \Lambda$ as $\varepsilon_h  \to 0$ and 
$ u_{\varepsilon_h} ( x_h ) \geq c > 0$. Therefore
\[ \min_{\partial \Lambda} \S \leq \S ( \bar{x} ) = \lim_{h \rightarrow
\infty }
   \S ( x_h ) = \min_{\Lambda} \S, \]
which contradicts our assumptions on $\Lambda$. Here we have used
the
continuity of $\S$, already proved. Hence we should just prove the
Claim 1. By
compactness, we assume without loss of generality that $x_h \rightarrow
\hat{x} \in \Lambda$. We proceed by contradiction, assuming therefore
that
\[ \S ( \hat{x} ) > \min_{\Lambda} \S . \]
Let $v_h ( x ) = u_{\varepsilon_h} ( x_h + \varepsilon_h x )$. By
elliptic
regularity, $\{ v_h \}$ converges strongly in $H^1 ( K )$ towards some
$v$,
for each compact set $K \subset \mathbb{R}^N$. Let $\chi$ be the weak*
limit (in $L^{\infty}$) of the sequence
$\{ \chi_{\Lambda} ( x_h + \varepsilon_h \cdot ) \}_{h \in
\mathbb{N} }$.
Clearly, $0 \leq \chi \leq 1 $. Therefore the function $v$ is a weak
solution
of the equation
\[ - \tmop{div} \left( J ( \hat{x} ) \nabla v \right) + V ( \hat{x} ) v
= g_0
   ( \cdot, v ) \quad \mbox{in } \mathbb{R}^N, \]
where
\[ g_0 ( x, s ) = \chi ( x ) f ( s ) + \left( 1 - \chi ( x ) \right) 
   \tilde{f} ( s ) . \]
Let now $E^h : H_V ( \Omega_h ) \rightarrow \mathbb{R}$ be the
functional
\[ E^h ( v ) = \frac{1}{2} \int_{\Omega_h} \left\langle J ( x_h +
   \varepsilon_h x ) \nabla v | \nabla v \right\rangle + \frac{1}{2}
   \int_{\Omega_h} V ( x_h + \varepsilon_h x ) v^2 - \int_{\Omega_h} G (
x_h +
   \varepsilon_h x, v ), \]
where $\Omega_h = \varepsilon_h^{- 1} \left( \Omega - x_h \right)$. 
Let us observe explicitly that 
\begin{equation}\label{eq:ugua}
E^h(v_h)=\eps^{-N}_h E^{\eps_h}(u_{\eps_h}).
\end{equation}
It is clear that $v_h$ is a critical point of $E^h$ in $H_V ( \Omega_h )$.

As remarked in \cite{dPF3} (see also \cite{S}), $v$ actually satisfies 
\begin{equation}\label{eq:v}
 - \tmop{div} \left( J ( \hat{x} ) \nabla v \right) + V ( \hat{x} ) v = f (   v ) \quad \textrm{ in } \RN. 
\end{equation}
Indeed, without loss of generality, we may suppose that $\chi ( x ) = \chi_{\{
x_1 < 0
\}} ( x )$ for all $x$. Consider the equation satisfied by $v$:
\[ 
- \tmop{div} \left( J ( \hat{x} ) \nabla v \right) + V ( \hat{x} ) v
=
   \chi_{\{ x_1 < 0 \}} ( x ) f ( v ) + \chi_{\{ x_1 > 0 \}} ( x )
\tilde{f} (   v ) \quad \textrm{ in } \RN, 
\]
and multiply it by $T_k \partial_{x_1} v$, where $T_k$ is a sequence of
smooth
functions such that $T_k ( x ) = 1$ if $| x | \leq k $, $T_k ( x ) = 0$
if $|
x | \geq 2 k $, and $| \nabla T_k ( x ) | = O ( 1 / k )$. 
Hence we get
\[ 
     \int_{\mathbb{R}^N} \left\langle J ( \hat{x} ) \nabla v | \nabla (
T_k
     \partial_{x_1} v ) \right\rangle + \int_{\mathbb{R}^N} V ( \hat{x}
) v
     T_k \partial_{x_1} v = \int_{\mathbb{R}^N} \varphi T_k
\partial_{x_1} v, 
\]
where we have set $\varphi ( x, v ) = \chi_{\{ x_1 < 0 \}} ( x ) f ( v (
x ))
+ \chi_{\{ x_1 > 0 \}} ( x ) \tilde{f} ( v ( x ))$. We observe that $v
\partial_{x_1} v = \frac{1}{2} \partial_{x_1} v^2$, so that after an
integration by parts we have
\begin{eqnarray*}
  \int_{\mathbb{R}^N} \left\langle J ( \hat{x} ) \nabla v |
\partial_{x_1} v
  \nabla T_k \right\rangle + \int_{\mathbb{R}^N} \left\langle J (
\hat{x} )
  \nabla v | T_k \partial_{x_1} \nabla v \right\rangle -
\int_{\mathbb{R}^N} V
  ( \hat{x} ) \frac{1}{2} v^2 \partial_{x_1} T_k &  & \\
  = \int_{\mathbb{R}^N} \varphi T_k \partial_{x_1} v . &  & 
\end{eqnarray*}
By a second integration by parts on the term $\int_{\mathbb{R}^N}
\left\langle
J ( \hat{x} ) \nabla v | T_k \partial_{x_1} \nabla v \right\rangle$, and
by
the Dominated Convergence Theorem, we pass to the limit as $k
\rightarrow +
\infty $ to get
\[ \int_{\mathbb{R}^N} \varphi ( x, v ) \partial_{x_1} v = 0 . \]
This implies easily that $v ( 0, x_2, \ldots, x_N ) \leq \ell $ for all
$(x_2, \ldots, x_N ) \in \mathbb{R}^{N - 1}$. It is now easy to check that we
can
choose
\[ \max \left\{ v ( \cdot ) - \ell, 0  \right\} \]
as a barrier for the equation satisfied by $v$, and thus prove that $v
\leq
\ell $ in $\mathbb{R}^N$. In particular, $\varphi ( x, v ( x )) = f ( v
( x
))$, so that $v$ is a solution of \eqref{eq:v}.
\vspace{0.5 cm}

\noindent\textit{Claim 2: there results
\begin{equation}\label{eq:semicont}
I_{\hat{x}} (v) \leq \liminf_{h\to\infty} E^h (v_h).
\end{equation}
}

Indeed, consider now the function
\[ \xi_h ( x ) = \frac{1}{2} \left\langle J ( x_h + \varepsilon_h x )
\nabla
   v_h | \nabla v_h \right\rangle + \frac{1}{2} V ( x_h + \varepsilon_h
x )
   v_h^2 = G ( x_h + \varepsilon_h x, v_h ) . \]
We already know that $v_h \rightarrow v$ strongly on compact sets. Therefore
\[ \lim_{h \rightarrow \infty } \int_{B_R} \xi_h = \frac{1}{2}
\int_{B_R}
   \left( \left\langle J ( \hat{x} ) \nabla v | \nabla v \right\rangle +
V (
   \hat{x} ) v^2 \right) - \int_{B_R} F ( v ) \]
for all $R > 0$. But $v \in H^1 ( \mathbb{R}^N )$, so that
\[ I_{\hat{x}} ( v ) - \frac{1}{2} \int_{B_R} \left( \left\langle J (
\hat{x}
   ) \nabla v | \nabla v \right\rangle + V ( \hat{x} ) v^2 \right) +
   \int_{B_R} F ( v ) = o ( 1 ) \]
as $R \rightarrow + \infty $. To show that $I_{\hat{x}} ( v ) \leq
\liminf_{h
\rightarrow \infty } E^h ( v_h )$ it is enough to prove that for all
$\delta >
0$ there exists $R > 0$ such that
\[ \liminf_{h \rightarrow \infty } \int_{\Omega_h \setminus B_R} \xi_h
\geq -
   \delta .  \]
We introduce again a cutoff function $\eta_R \in C^{\infty} (
\mathbb{R}^N )$
such that
\[ \eta_R = 0 \quad \tmop{in} B_{R - 1}, \]
\[ \eta_R = 1 \quad \tmop{in} \mathbb{R}^N \setminus B_R, \]
\[ | \nabla \eta_R | \leq C .  \]
We test the equation satisfied by $v_h$ against $\eta_R v_h$. After some
computations we get
\begin{gather*}
\liminf_{h\to\infty} \int_{\Omega_h\setminus B_R}\xi_h \geq - \frac{1}{2}\limsup_{h\to\infty} \left[ \int_{B_R\setminus B_{R-1}} \langle
J(x_h+\ge_h x)\nabla v_h \mid \nabla (\eta_R v_h)\rangle + \right.\\
\left. \int_{B_R\setminus B_{R-1}} V(x_h+\ge_h x)v_h^2 \eta_R - \int_{B_R\setminus B_{R-1}}g(x_h+\ge_h x,v_h)\eta_R v_h \right]=o(1)
\end{gather*}
as $R \to +\infty$.
This finally proves that $I_{\hat{x}} ( v ) \leq \liminf_{h
\rightarrow \infty } E^h ( v_h )$ and so Claim 2 holds.

We now complete the first part of the proof. 
First of all, from Lemma \ref{le:dis}, 
\eqref{eq:ugua} and from \eqref{eq:semicont}, it follows that
\[
I_{\hat{x}}(v) \leq \bar c=\inf_{\gamma\in\mathcal{P}_0} \sup_{0\leq t
\leq 1} I_0(\gamma (t)).
\]
On the other hand, since $v$ is a critical point of $I_{\hat{x}}$, by Proposition \ref{prop:ugua}, 
we have
\begin{eqnarray*}
I_{\hat{x}}(v) &\geq& \inf_{\gamma\in\mathcal{P}_{\hat{x}}} \sup_{0\leq
t \leq 1} I_{\hat{x}}(\gamma (t))\\
&=& \inf_{u\in \mathcal{N}_{\hat{x}}} I_{\hat{x}}(u) = \S (\hat{x}) \\
&>& \inf_{u\in \mathcal{N}_{z_0}} I_0 (u) = \S (z_0)\\
&=& \inf_{\gamma\in\mathcal{P}_0} \sup_{0\leq t \leq 1} I_0(\gamma
(t))=\bar c.
\end{eqnarray*}
This contradiction proves Claim 1 and so also the first part of the proof.

As regards the last statements of the proposition, these follow 
easily from the corresponding properties of solutions in \cite{dPF}. 
We just sketch the ideas. Let $\bar{z}=\lim_{\eps \to 0}x_\eps$ and 
take any critical point $u$ of $I_{\bar{z}}$.
By the change of variables $x \mapsto Tx$ introduced in the proof of Proposition 
\ref{prop:1}, if $v$ is a solution of equation \eqref{eq:rescaled}, namely
\[
-\varDelta v + V(\bar{z})v=f(v) \quad \textrm{ in } \RN,
\]
then $u(x)=v(Tx)$. It is well known by \cite{GNN} that solutions of \eqref{eq:rescaled} are radially 
symmetric and decreasing. In particular, $x=0$ is a nondegenerate maximum point of $u$. 
We are in a position to apply a reasoning similar to that of \cite{dPF}, page 133.
\end{proof}

We can now prove Theorem \ref{thm1}.

\bigskip

\begin{proof1}
By Proposition \ref{prop:key}, we know that if $\ge$ is small enough,
then
\[
u_\ge (x) < \ell \quad\hbox{for all $x\in\de\Lambda$}.
\]
The function $\left( u_\ge (\cdot)-\ell \right)^+ = \max\{u_\ge
(\cdot)-\ell,0\}$ belongs to $H_0^1(\Omega)$, so that we can test
the equation
\[
-\ge^2 \dv (J \nabla u_\ge) + V u_\ge = g(\cdot,u_\ge) \quad \textrm { in } \O
\]
against it. By the divergence theorem,
\begin{multline}\label{eq:final}
\ge^2 \int_{\Omega\setminus\Lambda} \langle J\nabla \left( u_\ge -\ell
\right)^+ \mid \nabla \left( u_\ge -\ell \right)^+ \rangle \\
+ \int_{\Omega\setminus\Lambda} \ell \Phi_\ge \left( u_\ge -\ell
\right)^+ + \int_{\Omega\setminus\Lambda} \Phi_\ge \left(\left( u_\ge -\ell
\right)^+\right)^2 =0,
\end{multline}
where we have set
\[
\Phi_\ge (x) = V(x) - \frac{g(x,u_\ge(x))}{u_\ge (x)}.
\]
The properties of $g$ imply that $\Phi_\ge > 0$ in $\Omega
\setminus\Lambda$. Therefore all the terms in \eqref{eq:final} must
vanish, and in particular
\[
u_\ge \leq \ell \quad \mbox{in $\Omega\setminus\Lambda$}.
\]
We conclude that $u_\ge$, for $\ge$ small enough, is actually a critical
point of $I^\eps$, and hence a solution of \eqref{eq:E}. The regularity
of $u_\ge$ follows again from \cite{G}. The last statement of the theorem follows immediately 
by Proposition \ref{prop:key}.
\end{proof1}

\begin{proof2}
First of all, we remark that, for any $u\in H^1(\RN)\setminus \{0\}$ and any fixed $z\in \partial \Lambda$, there results
\[
I_z (u) > I_{z_0}(u),
\]
because of \eqref{1} and \eqref{2}. We claim that there exists a ground state $v_z$ such that
\begin{equation}\label{eq:v_z}
I_z(v_z)=\inf_{u\in\mathcal{N}_z} I_z(u).
\end{equation}
Indeed, the usual change of variables $x \mapsto Tx$ rescales the functional $I_z$ by a constant factor $|\det T| > 0$. This reduces the search of a ground state for $I_z$ to the search of a ground state for the equation
\[
-\varDelta v + V(z)v = f(v) 
\mbox{\quad in $\RN$},
\]
whose existence follows easily from the results contained in \cite{BL} and \cite{ChL}. 
This proves the claim.

Let $z\in \de \L$ and $v_z$ as in \eqref{eq:v_z}. By Lemma \ref{le:retta} 
we know that there exists a positive constant $\tau$ such that $\tau v_z \in \Ne_{z_0}$. 
By our assumptions and \eqref{eq:max}, we easily get
\[
\S (z)= I_z(v_z) \geq I_z(\tau v_z) > I_{z_0}(\tau v_z) =\max_{t >0} I_{z_0} \left(t (\tau v_z)\right)
 \geq  \S(z_0).
\]
Therefore $\S(z_0) < \S(z)$ and so \eqref{eq:min} holds.
\end{proof2}

\section{Necessary condition for concentration}\label{sec:nec}

In this section we want to show that the function $\S$ also plays a
{\it necessary} r\^{o}le for the existence of concentrating solutions
of \eqref{eq:E}. We will give also a more general version of Theorem \ref{th:2}.

We suppose that $\Omega = \mathbb{R}^N$. Indeed, if $\Omega$ has a boundary, we do not expect that solutions must concentrate at critical points of $\Sigma$, but rather on critical point of some function connected to the geometry of $\partial\Omega$, see for example \cite{NW}.

\begin{theorem}\label{thm2}
Assume, in addition to assumptions {\bf (V)}, {\bf (J)},
{\bf (f1-4)}, that $V$ is bounded from above, and $\Omega=\mathbb{R}^N$. Let $\{u_{\ge_j}\}$ be a
sequence of solutions of \eqref{eq:E} such that for all $\ge >0$ there
exist $\rho >0$ and $j_0>0$ such that for all $j \geq j_0$ and for all points
$x$ with $|x-z_0| \geq \ge_j \rho$, there results
\[
u_{\ge_j}(x) \leq \ge.
\]
If, for all $z\in \RN$, the functional $I_z$ has only one positive ground-state 
(up to translations), then $z_0$ is a critical point of $\S$.
\end{theorem}

Before proving the theorem, we recall a recent version of Pucci--Serrin variational identity 
for lipschitz continuous solutions of a general class of Euler equations (see \cite{DMS}).

\begin{theorem}\label{th:PS}
Let ${\cal L} \colon \RN\times \R \times \RN \to \R$ be a $C^1$ function such that 
the function $\xi \mapsto {\cal L}(x,s,\xi)$ is strictly convex for every $(x,s)\in \RN \times \R$.
Let $\varphi \in L^\infty_{{\rm loc}}(\RN)$.
Let $u\colon \RN \to \R$ be a locally Lipschitz weak solution of 
\[
-\dv \left( \de_\xi \cl(x,u,\nabla u)\right) + \de_s \cl(x,u,\nabla u)=\varphi 
\quad {\rm in } \;\RN. 
\]
Then
\begin{multline*}
\sum\limits^N_{i,\,k=1}\int_\RN \de_i h^k \de_{\xi_i}\cl(x,u,\nabla u) \de_k u
\\
- \int_\RN \big[(\dv h)\cl(x,u,\nabla u) + h \cdot \de_x  \cl(x,u,\nabla u)     \big]=
\int_\RN (h\cdot \nabla u )\varphi,
\end{multline*}
for all $h \in C^1_{\rm c}\left(\RN, \RN \right)$.
\end{theorem}

\begin{proofnec}
To save notation, we write $u_j$ instead of $u_{\ge_j}$. Define
$w_j(x)=u_j(z_0+\ge_j x)$. Therefore
\begin{equation}\label{eq:wj}
- \dv \left(J(z_0+\ge_jx)\nabla w_j\right) + V(z_0+\ge_j x)w_j - f(w_j)=0.
\end{equation}
By assumption, $w_j$ decays to zero uniformly with respect to
$j\in\mathbb{N}$. It is not difficult to build an exponential barrier
for $w_j$, proving in this way that $w_j$ decays to zero exponentially
fast at infinity. By elliptic regularity (see \cite{G}), the sequence
$\{w_j\}$ converges in $C^2_{\rm loc}$ to a solution $w_0$ of the
equation
\[
-\dv \left(J(z_0) \nabla w_0 \right) + V(z_0) w_0 -f(w_0)=0.
\]

Let us apply Theorem \ref{th:PS} to \eqref{eq:wj}, with
\begin{eqnarray*}
\cl(x,s,\xi) & = & \frac{1}{2}\langle J(z_0+\ge_jx) \xi \mid \xi \rangle+ 
\frac{1}{2} V(z_0+\ge_jx) s^2 - F(s),
\\
h(x) & = & (T(\eps x), \,0, \ldots, \,0),
\\
\varphi (x) & = & 0,
\end{eqnarray*}
where $T \in C^1_{{\rm c}} (\RN)$ such that $T(x)=1$ if $|x|\leq 1$ and $T(x)=0$ if $|x|\geq 2$. 

By Theorem \ref{th:PS}, we have:\footnote{We denote by $J_i (x)$ the $i$--th row of $J(x)$.}
\begin{multline*}
\sum\limits^N_{i=1}\int_\RN \eps \de_i T(\eps x)\langle J_i(z_0+\ge_jx) \mid \nabla w_j \rangle \de_1 w_j
\\
- \int_\RN \eps \de_1 T(\eps x)\left[\frac{1}{2}\langle J(z_0+\ge_jx) \nabla w_j \mid \nabla w_j \rangle
+ \frac{1}{2} V(z_0+\ge_jx)w_j^2 -F(w_j) \right]
\\
- \int_\RN \eps T(\eps x)\left[\frac{1}{2}\langle \de_1J(z_0+\ge_jx) \nabla w_j \mid \nabla w_j \rangle
+ \frac{1}{2} \de_1 V(z_0+\ge_jx)w_j^2 \right]=0
\end{multline*}

Passing to the limit in the previous relations, as $\eps \to 0$, we get:
\[
\frac{1}{2}\int_\RN \left(\langle \de_1 J(z_0)\nabla w_0 \mid \nabla
w_0 \rangle + \de_1 V(z_0) |w_0|^2
\right)=0.
\]
The proof is complete once we recall that if $S^{z_0}$ consists of just one
element, then, by Proposition \ref{prop:1}, $\de_1 \S(z_0)=0$. The proof for the other 
partial derivatives is identical. 
\end{proofnec}

\begin{remark}
Let us observe that Theorem \ref{th:2} is an immediate consequence of Theorem \ref{thm2}. Indeed 
by \cite{ChL} we know that the problem \eqref{eq:E1} with frozen coefficients has only one 
positive ground-state (up to translations).
\end{remark}

\section{Existence via perturbation method}\label{sec:ex-per}

We have seen in the previous sections that the penalization technique of del~Pino and Felmer 
provides at least a solution of \eqref{eq:E} if the auxiliary map $\Sigma$ possesses a minimum. 
Moreover, we could also treat the case of maximum point of $\S$ under some more restrictive, 
global, assumptions on the potentials $J$ and $V$. In the present section we show that for 
\eqref{eq:E1} it is possible to find at least a solution just by differential methods if 
there exists a local maximum or minimum of $\S$. 
More precisely, we will apply the perturbation technique in critical point theory as 
developed in \cite{AMS}. Since this approach deals with the local behavior of the 
potentials $J$ and $V$, we need a better knowledge about the derivatives of the potentials.

In addition to hypotheses {\bf (V)} and {\bf (J)}, in this section we will always assume: 
\begin{description}
\item[(V1)] $V\in C^{2}(\RN, \R)$, $V$ and $D^{2}V$ are bounded;
\item[(J1)] $J\in C^{2}(\RN,\R^{N\times N})$, $J$ and $D^{2}J$ are bounded.
\end{description}
 
Since we follow closely \cite{AMS}, 
we will skip some proof and we will give only the sketch of some others.

Without loss of generality we can assume that $V(0)=1$. Moreover, using the change of variables 
introduced in the proof of Proposition \ref{prop:1}, we can assume also $J(0)=I$, where $I$ is the 
identity matrix of order $N\times N$.

Performing the change of variable $x\mapsto \eps x$, equation 
\[
-\ge^2\dv \big( J(x)\nabla u \big) + V(x)u=u^p \quad \textrm{ in } \RN
\]
becomes
\begin{equation}\label{eq:P}
- \dv\left(J(\eps x)\nabla u\right)+V(\eps x)u=u^p  \quad \textrm{ in } \RN.
\end{equation}
Solutions of (\ref{eq:P}) are the critical points $u\in H^1(\RN)$ of
\[
f_\eps (u)=f_0(u)+ 
\frac{1}{2}\int_{\RN}\langle \left( J(\eps x)- I) \right) \nabla u \mid \nabla u \rangle dx+
\frac{1}{2}\int_{\RN} \left( V(\eps x) -1 \right)u^2dx,
\]
where
\[
f_0(u)=\frac{1}{2}\|u\|^2 - \frac{1}{p+1}\int_{\RN}u^{p+1}dx,
\]
and $\|u\|^2= \int_{\RN} |\nabla u|^2 +u^2$.
The solutions of (\ref{eq:P}) will be found near a solution of
\begin{equation}\label{eq:xi}
- \dv\left(J(\eps \xi)\nabla u\right) +V(\eps \xi)u=u^p,
\end{equation}
for an appropriate choice of $\xi\in \RN$. 

The solutions of (\ref{eq:xi}) are critical points of the functional
\begin{equation}\label{eq:F}
F^{\eps\xi}(u)=f_{0}(u)+
\frac{1}{2}\int_{\RN}\langle \left( J(\eps \xi)- I \right) \nabla u \mid \nabla u \rangle +
\frac{1}{2}\left(V(\eps \xi)-1 \right)\int_{\RN}u^2dx
\end{equation}
and can be found explicitly. First of all, by the usual change of variables, $x \mapsto T(\eps \xi)$, equation \eqref{eq:P} becomes
\[
- \varDelta v +V(\eps \xi)v=v^p.
\]
Let $U$ denote the unique, positive,
radial solution of
\begin{equation}\label{eq:unp}
 - \varDelta u + u =u^p ,\qquad u \in H^1(\RN).
\end{equation}
Then a straight calculation shows that $\a U(\b T x)$ solves (\ref{eq:P})
whenever
\[
\a = \a(\eps\xi)=[V(\eps\xi)]^{1/(p-1)}, \quad 
\b=\b(\eps\xi)= [V(\eps\xi)]^{1/2}\quad {\rm and}\quad T=T(\eps \xi).
\]
We set
\begin{equation}\label{eq:zU}
z^{\eps\xi}(x)=\a(\eps\xi)U\big(\b(\eps\xi)T(\eps \xi) x\big)
\end{equation}
and
\[
Z^{\eps}=\{z^{\eps\xi}(x-\xi):\xi\in \RN\}.
\]
When there is no possible misunderstanding, we will write $z$, resp. $Z$, 
instead of $z^{\eps\xi}$, resp $Z^{\eps}$. 
We will also use the notation $z_{\xi}$ to denote the function 
$z_{\xi}(x):=z^{\eps\xi}(x-\xi)$. 
Obviously all the 
functions in $z_{\xi}\in Z$ are solutions of (\ref{eq:xi}) or, equivalently, 
critical points of $F^{\eps\xi}$. 

\begin{remark}
Before we proceed, a remark is in order concerning the definition of the manifold $Z^\ge$. Indeed, there is a lot of freedom in the choice of the diagonalizing matrix $T$. Moreover, $Z^\ge$ should be a regular manifold. We claim that, thanks to the uniform ellipticity assumption on $J$, see {\bf (J)}, it is possible to choose $T(\ge\xi)$ with the same regularity as $J$ itself. We will not supply a complete proof of this fact. However, the best way to convince oneselves of this is to remember the celebrated Householder algorithm that diagonalises a symmetric matrix $J$ by means of arithmetic operations on the rows and the columns of $J$. We refer to \cite{Str} for an explanation of the method. Each of these operations corresponds to an orthogonal change of variables which preserves the uniform ellipticity of $J$, and at each step the only possible lack of regularity can be due to the division by an entry on the main diagonal of $J$. By {\bf (J)}, each such entry is a function of $\ge\xi$ strictly bounded away from zero, so that it cannot introduce any singularity in the algorithm. A repeated application of this argument can now applied to prove the regularity of $T$.
\end{remark}

For future references let us point out some estimates. First of all, by straightforward 
calculations, we get:
\begin{equation}\label{eq:partialz}
\partial_{\xi}z^{\eps\xi}(x-\xi)=-\partial_{x}z^{\eps\xi}(x-\xi)+O(\eps).
\end{equation}
Moreover, using {\bf (J1)} and {\bf (V1)}, we can infer that $\nabla f_{\eps}(z_\xi)$ 
is close to zero when $\eps$ is small. Indeed we have:
\begin{equation}\label{eq:1}
\|\nabla f_{\eps}(z_{\xi})\|\leq C
\left(\eps |D J(\eps\xi)| +
\eps |\nabla V(\eps\xi)|+\eps^{2}\right),\quad C>0.
\end{equation}

In the next lemma we will show that $D^{2}f_{\eps}$ is invertible on
$\left(T_{z_\xi}Z^\eps \right)^{\perp}$, where $T_{z_\xi} Z^\eps$ 
denotes the tangent space to $Z^\eps$ at $z_\xi$.

Let $L_{\eps,\xi}:(T_{z_{\xi}}Z^{\eps})^{\perp}\to
(T_{z_{\xi}}Z^{\eps})^{\perp}$ denote the operator defined by setting
$(L_{\eps,\xi}v|w)= D^{2}f_{\eps}(z_{\xi})[v,w]$.

\begin{lemma}\label{lem:inv}
Given $\overline{\xi}>0$ there exists $C>0$ such that for $\eps$ small enough
one has that
\begin{equation}\label{eq:inv}
|(L_{\eps,\xi}v|v)|\geq C \|v\|^{2},\qquad \forall\;|\xi|\leq
\overline{\xi},\;\forall\; v\in(T_{z_{\xi}}Z^{\eps})^{\perp}.
\end{equation}
\end{lemma}

\begin{proof}
We recall that $T_{z_\xi} Z^\eps = {\rm span} \{\de_{\xi_1}z_\xi, \ldots, \de_{\xi_N}z_\xi \}$. 
Let ${\cal V}= {\rm span} \{z_\xi,$ $\de_{x_1}z_\xi, \ldots, \de_{x_N}z_\xi \}$, by (\ref{eq:partialz}) 
it suffices to prove (\ref{eq:inv}) for all $v\in {\rm span}\{z_{\xi},\phi\}$, where $\phi$ is 
orthogonal to ${\cal V}$. Precisely we shall prove that there exist $C_{1},C_{2}>0$ such that
for all $\eps>0$ small and all $|\xi|\leq \overline{\xi}$
one  has:
\begin{eqnarray}
(L_{\eps,\xi}z_{\xi}|z_{\xi})& \leq & - C_{1}< 0.      \label{eq:neg} 
\\
(L_{\eps,\xi}\phi|\phi)&\geq & C_{2} \|\phi\|^2.        \label{eq:claim}
\end{eqnarray}

The proof of (\ref{eq:neg}) follows easily from the fact that $z_{\xi}$ is a Mountain Pass 
critical point of $F^{\eps \xi}$ and so from the fact that, given $\overline{\xi}$, there exists $c_0>0$ 
such that for all $\eps>0$ small
and all $|\xi|\leq \overline{\xi}$ one finds:
\begin{equation*}
D^2 F^{\eps\xi}(z_{\xi})[z_{\xi},z_{\xi}] < -c_0< 0.
\end{equation*}

Let us prove (\ref{eq:claim}). As before, the fact that
$z_{\xi}$ is a Mountain Pass critical point of $F^{\eps \xi}$ implies that
\begin{equation}\label{eq:15}
D^2 F^{\eps\xi}(z_{\xi})[\phi,\phi]>c_1 \|\phi\|^2 \quad \forall \phi \perp {\cal V}.
\end{equation}
Let $R\gg 1$ and consider a radial smooth function
$\chi_{1}:\RN\mapsto \R$ such that
\begin{equation*}
\chi_{1}(x) = 1, \quad \hbox{ for } |x| \leq R; \qquad
\chi_{1}(x) = 0, \quad \hbox{ for } |x| \geq 2 R;
\end{equation*}
\begin{equation*}
|\nabla \chi_{1}(x)| \leq \frac{2}{R}, \quad \hbox{ for } R \leq |x| \leq 2 R.
\end{equation*}
We also set $ \chi_{2}(x)=1-\chi_{1}(x)$.
Given $\phi$ let us consider the functions
\[
\phi_{i}(x)=\chi_{i}(x-\xi)\phi(x),\quad i=1,2.
\]
Therefore we need to evaluate the three terms in the
equation below:
\[
(L_{\eps,\xi}\phi|\phi)=
(L_{\eps,\xi}\phi_{1}|\phi_{1})+
(L_{\eps,\xi}\phi_{2}|\phi_{2})+
2(L_{\eps,\xi}\phi_{1}|\phi_{2}).
\]
Using \eqref{eq:15} and the definition of $\chi_i$, we easily get
\begin{eqnarray*}
(L_{\eps,\xi}\phi_{1}|\phi_{1}) & \geq & \ge c_{1}\|\phi_{1}\|^{2}-
\eps c_{2}\|\phi\|^{2}+o_{R}(1)\|\phi\|^{2},
\\
(L_{\eps,\xi}\phi_{2}|\phi_{2}) & \geq & c_{3} \|\phi_{2}\|^{2}+o_{R}(1)\|\phi\|^{2},
\\
(L_{\eps,\xi}\phi_{1}|\phi_{2}) & \geq & o_{R}(1)\|\phi\|^{2}.
\end{eqnarray*} 
Therefore, since
\[
\| \phi \|^2 = \| \phi_1 \|^2 + \| \phi_2 \|^2 + 2\int_\RN \chi_{1}\chi_{2}(\phi^{2}+|\nabla \phi|^{2}) 
+ o_R(1)\| \phi \|^2,
\]
we get
\[
(L_{\eps,\xi}\phi|\phi)\geq c_{4}\|\phi\|^{2}-c_{5} R \eps\|\phi\|^{2}+o_{R}(1)\|\phi\|^{2}.
\]
Taking $R = \ge^{-1/2}$, and choosing $\ge$ small, (\ref{eq:claim}) follows.

%
This completes the proof of the lemma.
\end{proof}

We will show that the existence of critical points of $f_{\eps}$ can be reduced to 
the search of critical points of an auxiliary finite dimensional functional. First of all we will make 
a Liapunov-Schmidt reduction, and successively we will study the behavior of an 
auxiliary finite dimensional functional.

\begin{lemma}\label{lem:w}
For $\eps>0$ small and $|\xi|\leq \overline{\xi}$ there exists a unique
$w=w(\eps,\xi)\in
(T_{z_\xi} Z)^{\perp}$ such that
$\nabla f_\eps (z_\xi + w)\in T_{z_\xi} Z$.
Such a $w(\eps,\xi)$ is of class $C^{2}$, resp.  $C^{1,p-1}$, with respect to $\xi$, provided that
 $p\geq 2$, resp. $1<p<2$.
Moreover, the functional $\Phi_\eps (\xi)=f_\eps (z_\xi+w(\eps,\xi))$ has
the same regularity of $w$ and satisfies:
 $$
\nabla \Phi_\eps(\xi_0)=0\quad \Longleftrightarrow\quad \nabla
f_\eps\left(z_{\xi_0}+w(\eps,\xi_0)\right)=0.
$$
\end{lemma}

\begin{proof}
Let $P=P_{\eps\xi}$ denote the projection onto $(T_{z_\xi} Z)^\perp$. We want
to find a solution $w\in (T_{z_\xi} Z)^{\perp}$ of the equation
$P\nabla f_\eps(z_\xi +w)=0$.  One has that $\nabla f_\eps(z+w)=
\nabla f_\eps (z)+D^2 f_\eps(z)[w]+R(z,w)$ with $\|R(z,w)\|=o(\|w\|)$, uniformly
with respect to $z=z_{\xi}$, for $|\xi|\leq \overline{\xi}$. Therefore, our equation is:
\[
L_{\eps,\xi}w + P\nabla f_\eps (z)+PR(z,w)=0.
\]
According to Lemma \ref{lem:inv}, this is equivalent to
\[
w = N_{\eps,\xi}(w), \quad \mbox{where}\quad
N_{\eps,\xi}(w)=-L_{\eps,\xi}\left( P\nabla f_\eps (z)+PR(z,w)\right).
\]
By \eqref{eq:1} it follows that
\begin{equation}\label{eq:N}
\|N_{\eps,\xi}(w)\|\leq c_1 \left(\eps |D J(\eps\xi)| +\eps|\nabla V(\eps\xi)|+\eps^2\right)+ o(\|w\|).
\end{equation}
Then one readily checks that $N_{\eps,\xi}$ is a contraction on some ball in
$(T_{z_\xi} Z)^{\perp}$
provided that $\eps>0$ is small enough and $|\xi|\leq \overline{\xi}$.
Then there exists a unique $w$ such that $w=N_{\eps,\xi}(w)$.  Let us
point out that we cannot use the Implicit Function Theorem to find
$w(\eps,\xi)$, because the map $(\eps,u)\mapsto P\nabla f_\eps (u)$ fails to be
$C^2$.  However, fixed $\eps>0$ small, we can apply the Implicit
Function Theorem to the map $(\xi,w)\mapsto P\nabla f_\eps (z_\xi + w)$.
Then, in particular, the function $w(\eps,\xi)$ turns out to be of class
$C^1$ with respect to $\xi$.  Finally, it is a standard argument, see
\cite{AB,ABC}, to check that the critical points of $\Phi_\eps
(\xi)=f_\eps (z+w)$ give rise to critical points of $f_\eps$.
\end{proof}

\begin{remark}\label{rem:w}
From (\ref{eq:N}) it immediately follows that:
\begin{equation}\label{eq:w}
\|w\|\leq C \left(\eps |D J(\eps\xi)| +\eps |\nabla V(\eps\xi)|+\eps^2\right),
\end{equation}
where $C>0$.
\end{remark}

With easy calculations (see Lemma 4 of \cite{AMS}), we can give an estimate of 
the derivative $\partial_\xi w$.

\begin{lemma}\label{lem:Dw}
One has that:
\begin{equation}\label{eq:Dw}
\|\partial_\xi w\|\leq c \left(\eps |D J(\eps\xi)|+ 
\eps |\nabla V(\eps\xi)|+O(\eps^2)\right)^\gamma,
\end{equation}
with $c>0$ and $\gamma=\min\{1,p-1\}$.
\end{lemma}

Now we will use the estimates on $w$ and $\partial_\xi w$ established above to find
an expansion of $\nabla \Phi_\eps (\xi)$, where $\Phi_\eps (\xi)= f_\eps(z_\xi
+w(\eps,\xi))$. In the sequel, to be short, we will often write $z$ instead of $z_\xi$
and $w$ instead of $w(\eps,\xi)$. It is always understood that $\eps$ is taken in such a 
way that all the results discussed previously hold.

We have:
\begin{multline*}
\Phi_\eps (\xi) = \frac{1}{2}\|z+w\|^2 + 
\frac{1}{2}\int_{\RN}\langle \left(J(\eps x)-I\right)\nabla (z+w) \mid \nabla (z+w) \rangle +
\\
\frac{1}{2}\int_{\RN} \left(V(\eps x) -1\right) (z+w)^2
-\frac{1}{p+1}\int_{\RN} (z+w)^{p+1}.
\end{multline*}
Since $-\dv\left( J(\eps \xi) \nabla z\right)+V(\eps\xi)z=z^p$, we infer that
\begin{eqnarray*}
\|z\|^2 & = & -\int_{\RN}\langle \left(J(\eps \xi)-I\right)\nabla z \mid \nabla z \rangle  
-\left(V(\eps\xi) -1\right) \int_{\RN} z^2 + \int_{\RN} z^{p+1},
\\
(z|w) & = & -\int_{\RN}\langle \left(J(\eps \xi)-I\right)\nabla z \mid \nabla w \rangle 
-\left(V(\eps\xi)-1\right)\int_{\RN} zw + \int_{\RN} z^pw.
\end{eqnarray*}
Then we find:
\begin{multline*}
\Phi_\eps (\xi)=
\left(\frac{1}{2}-\frac{1}{p+1}\right) \int_{\RN} z^{p+1}+
\\
\frac{1}{2}\int_{\RN}\langle \left(J(\eps x) - J(\eps \xi)\right)\nabla z \mid \nabla z \rangle
+\frac{1}{2}\int_{\RN}\left[V(\eps x)-V(\eps\xi)\right]z^{2}+
\\
\int_{\RN}\langle \left(J(\eps x) - J(\eps \xi)\right)\nabla z \mid \nabla w \rangle
+\int_{\RN}\left[V(\eps x)-V(\eps\xi)\right]zw+
\\
\frac{1}{2}\int_{\RN}\langle J(\eps x)\nabla w \mid \nabla w \rangle
+\frac{1}{2}\int_{\RN} V(\eps x)w^2+ 
\\ 
\frac{1}{2} \|w\|^2 
- \frac{1}{p+1}\int_{\RN}
\left[(z+w)^{p+1}-z^{p+1}-(p+1)z^pw\right].
\end{multline*}
Since $z(x)=\a(\eps\xi)U\big(\b(\eps\xi)T(\eps \xi)x\big)$, see (\ref{eq:zU}),
it follows
\[
\int_{\RN} z^{p+1}dx= C_0 V(\eps\xi)^{\frac{p+1}{p-1}-\frac{N}{2}} 
\left(\det J(\eps \xi )\right)^{\frac{1}{2}} = C_0 \G(\eps \xi),
\]
where $C_0= \int_{\RN} U^{p+1}$ and $\G$ is the auxiliary function introduced in \eqref{eq:Gamma}.
Letting $C_1= C_0 [1/2 -1/(p+1)]$ and
recalling the estimates (\ref{eq:w}) and (\ref{eq:Dw}) on $w$ and $\nabla_\xi w$, respectively,
we readily find:
\begin{equation}\label{eq:Phi}
\Phi_\eps (\xi)= C_1 \G(\eps\xi) + \rho_\eps(\xi),
\end{equation}
where $|\rho_\eps(\xi)|\leq {\rm const} \left(\eps |D J (\eps\xi)|+
\eps |\nabla V(\eps\xi)| +\eps^2\right)$ and
\begin{equation}\label{eq:DF}
\nabla \Phi_\eps (\xi)= C_1 \eps\nabla \G (\eps\xi) + \eps^{1+\gamma} R_\eps(\xi),
\end{equation}
where $|R_\eps(\xi)|\leq {\rm const}$ and $\gamma=\min\{1,p-1\}$.\hfill\par

\begin{remark}
We highlight that, as observed in Remark \ref{rem:SG},
$C_1\G =\S$, where $\S$ is the ground-state function.
\end{remark}

Now we can prove Theorem \ref{th:ex} at least in the case $(a)$. The other 
is an easy consequence of Theorem \ref{th:main}.

\begin{proofex}
Let $z_0$ be a minimum point of $\G$ (the other case is similar) and let 
$\L \subset \RN$ be a compact neighborhood of $z_0$ such that 
\begin{equation*}
\min_\L \G<\min_{\de \L}\G.
\end{equation*}
By \eqref{eq:Phi}, it is easy to see that for $\eps$ sufficiently small, there results:
\[
\min_{\L} \Phi(\cdot/\eps)<\min_{\de \L}\Phi(\cdot/\eps).
\]
Hence, $\Phi(\cdot/\eps)$ possesses a critical point $\xi$ in $\L$. By Lemma \ref{lem:w} we have that 
$u_{\eps, \xi}=z^\xi(\cdot-\xi/\eps)+ w(\eps, \xi)$ is a critical point of $f_\eps$ and so a solution of 
problem \eqref{eq:P}. Therefore
\[
u_{\eps, \xi}(x/\eps)\simeq z^\xi\left(\frac{x-\xi}{\eps} \right)
\]
is a solution of (\ref{eq:E1}). This $\xi$ converges to some $\bar{\xi}$ as $\eps \to 0$, but by 
\eqref{eq:DF} it follows that $\bar{\xi}=z_0$.
\end{proofex}

\section{Existence of multiple solutions}\label{sec:mult}

In this section we will study the problem of the multiplicity of solutions. In the first subsection 
we will prove that under a more stringent assumption on the function $\Sigma$, 
our problem \eqref{eq:E} possesses generically more than one solution. In the second subsection, instead, 
we will deal with the problem \eqref{eq:E1} and, as done in Section \ref{sec:ex-per}, we 
will treat also the case of maximum points for $\S$.

\subsection{Using penalization method}
Since our arguments are inspired to those of \cite{CL} and \cite{dPF2}, we will skip some easy details.

Let $c_0=\min_{\RN} \Sigma (z)$. Let $M \subset \S^{-1} (c_0) \cap \O$.

We state our main result for multiple solutions.

\begin{theorem}\label{th:multi}
Suppose {\bf (V)}, {\bf (J)}, {\bf (f1-4)}. Suppose that $M$ is compact 
and let $\L \subset \O$ be the closure of a bounded neighborhood of $M$ such that 
$c_0 < \inf_{\partial \Lambda} \Sigma$.

Suppose, in addition, that there exists a point $z_0\in M$ such that:
\begin{description}
\item[(V2)] $V(z_0)=\min_\L V$;
\item[(J2)] the matrix $J(z)-J(z_0)$ is positive-definite for all $z\in\RN$.
\end{description}
Then there exists $\eps (\L)>0$ such that, for any $\eps<\eps (\L)$, 
problem \eqref{eq:E} has at least $\cat(M, \L)$ solutions 
concentrating at some points of $M$. Here $\cat(M, \L)$ denotes the 
Lusternik-Schnirelman category of $M$ with respect to $\L$.
\end{theorem}


The proof of theorem \ref{th:multi} requires some preliminary lemmas. The main ingredient is the following result in abstract critical point theory (see for example \cite{Wil}).

\begin{theorem}\label{th:cat}
Let $X$ be a complete Riemannian manifold of class $C^{1,1}$, and assume that $\phi\in C^1(X)$ is bounded from below. Let
\[
-\infty < \inf_X \phi <a<b<+\infty.
\]
Suppose that $\phi$ satisfies the Palais--Smale condition on the sublevel $\{u\in X \mid \phi (u)\leq b\}$ and that $a$ is not a critical level for $\phi$. Then the number of critical points of $\phi$ in $\phi^a=\{u\in X \mid \phi(u)\leq a\}$ is at least $\cat(\phi^a, \phi^a)$.
\end{theorem}

We shall apply this theorem to the penalized functional $E^{\ge}$, introduced in \eqref{eq:E-eps}, 
constrained to its Nehari manifold $\Ne^\ge$, so that it satisfies (PS) and it is bounded from below. The crucial step is therefore to link the topological richness of the sublevels of $E^\ge$ with that of $M$. For this purpose we make use of the following elementary result. For the proof we refer to \cite{BC}.

\begin{lemma}\label{le:cat}
Let $H$, $\Omega^+$, $\Omega^-$ be closed sets with $\Omega^- \subset \Omega^+$; let $\beta \colon H \to \Omega^+$, $\psi \colon \Omega^- \to H$ be two continuous maps such that $\beta \circ\psi$ is homotopically equivalent to the embedding $j\colon \Omega^- \to \Omega^+$. Then $\cat(H, H) \geq \cat(\Omega^-,\Omega^+)$.
\end{lemma}

Let $\eta >0$ be a smooth, non-increasing cut--off function, defined in $[0,+\infty)$, such that $\eta (|x|)=1$ if $x\in\Lambda$, and $|\eta'|\leq C$ for some $C>0$. For any $\xi\in M$ let
\[
\psi_{\ge,\xi} \colon x\mapsto \eta (|x-\xi|) \, \omega \left(\frac{x-\xi}{\ge}\right),
\]
where $\omega$ is a positive ground state of the functional $I_\xi$. 
Now define $\Phi_\ge \colon M \to \Ne^\ge$ by
\[
\Phi_\ge (\xi)=\theta_\ge \psi_{\ge,\xi},
\]
where $\theta_\ge\in\R$ is such that $\theta_\ge \psi_{\ge,\xi}\in \Ne^\ge$. 
By Lemma \ref{le:retta} with minor changes, we infer that there exists such a $\t_\eps$.

\begin{lemma}\label{le:c0}
Uniformly in $\xi\in M$ we have
\[
\lim_{\ge\to 0} \ge^{-N} E^{\ge} (\Phi_\ge (\xi))=c_0.
\]
\end{lemma}

\begin{proof}
The proof is similar to the one of Lemma 4.1 in \cite{CL}, taking into account the monotonicity property {\bf (f4)} of $f$ and the fact that $\xi\in M\subset \Lambda$.
\end{proof}

We now construct a second auxiliary map which proves to be useful for the comparison of the topologies of $M$ and of the sublevels of $E^\ge$.

Let $R>0$ be such that $\L \subset \{x\in\RN \colon |x| \leq R\}$. Let $\chi\colon \RN \to \RN$ be defined by
\begin{equation*}
\chi (x)=
\begin{cases}
x &\text{ for $|x| \leq R$} \\
\frac{R x}{|x|} &\text{ for $|x| > R$}.
\end{cases}
\end{equation*}
Finally, define $\beta \colon \Ne^\ge \to \RN$ by
\[
\beta (u) = \frac{\int_\RN \chi \cdot |u|^2}{\int_\RN |u|^2}.
\]
As in \cite{CL}, it is easy to show that $\beta (\Phi_\ge(\xi))=\xi +o(1)$ as $\ge \to 0$, uniformly with respect to $\xi\in M$.

Let us define a suitable sublevel of $E^\ge$:
\[
\tilde{\Ne}^\ge = \left\{ u\in \Ne^\ge \colon E^\ge (u) \leq \ge^N (c_0+o(1)) \right\}.
\]

As already stated, we know that $E^\ge$ verifies the (PS) condition at all levels.

\begin{lemma}\label{le:beta}
Let $\tilde{\L}$ a sufficiently small homotopically equivalent neighborhood of $\L$. 
For all $\ge$ sufficiently small, we get
\[
\beta ( \tilde{\Ne}^\ge ) \subset \tilde{\L}.
\]
\end{lemma}

\begin{proof}
The proof proceeds by contradiction. If the claim does not hold, then we may find sequences 
$\{\ge_n\}$, $\{u_n\}$ such that $\ge_n\to 0$, $u_n\in \tilde{\Ne}^{\ge_n}$ but 
$\beta(u_n) \notin \tilde{\L}$. We claim that
\begin{equation}\label{claim1}
\lim_{n\to\infty} \ge_n^{-N} \int\limits_{\Omega \setminus \tilde{\L}} |u_n|^2 =0.
\end{equation}
Indeed, since $u_n\in \Ne^{\ge_n}$, we have that
\[
E^{\ge_n} (u_n) \geq E^{\ge_n}(t u_n)
\]
for any $t>0$. Let us set
\[
\tilde{E}_n (v)=\frac{1}{2} \int_{\tilde{\L}} \ge_n^2 \langle J(x)\nabla v \mid \nabla v \rangle +
 V(x)|v|^2 - \int_{\tilde{\L}} G(x,v)\, dx.
\]
Choose $t_n >0$ such that
\[
\tilde{E}_n(t_n u_n)=\max_{t>0} \tilde{E}_n (t u_n).
\]
Since $u_n \in \tilde{\Ne}^{\ge_n}$ and that fact that
\[
\frac{V(x)}{2} u^2 - G(x,u) \geq C u^2
\]
for all $x\in \Omega \setminus \tilde{\L}$ and all $u>0$, we obtain
\begin{equation}\label{eq:1.22}
\tilde{E}_n(t_n u_n) + C t_n^2 \int\limits_{\Omega \setminus \tilde{\L}} |u_n|^2 \leq \ge_n^{N}(c_0 + o(1)).
\end{equation}
From our assumptions on $V$, $J$ and $f$, since $E^{\ge_n}(u_n) \leq C \ge_n^N$ and $u_n \in \Ne^{\ge_n}$, we see that
\begin{equation}\label{eq:1.24}
\int_\Omega \ge_n^2 |\nabla u_n|^2 + |u_n|^2 \leq C \ge_n^N.
\end{equation}
Set $v_n \colon x \mapsto t_n u_n (\ge_n x)$. From the definition of $t_n$ it follows
\begin{align*}
\int\limits_{\ge_n^{-1}\tilde{\L}} \langle J(\eps_n x )\nabla v_n \mid \nabla v_n \rangle + V(\ge_n x)|v_n|^2 &=
\int\limits_{\ge_n^{-1}\tilde{\L}} g(\ge_n x,v_n)v_n \\
&\leq \int\limits_{\ge_n^{-1}\tilde{\L}} C |v_n|^{p+1} + \rho |v_n|^2,
\end{align*}
where $\rho >0$ can be taken arbitrarily small. Now, Sobolev's theorem yields that
\[
\int\limits_{\ge_n^{-1}\tilde{\L}} |v_n|^{p+1} \leq 
C \Big( \,  \int\limits_{\ge_n^{-1}\tilde{\L}}|\nabla v_n|^2 +|v_n|^2 \Big)^{\frac{p+1}{2}}
\]
and the constant $C$ can be taken the same for all $n$, since it generally depends 
only on the geometry of the domain of integration but not on its volume. Combining 
the two last inequalities, since $J$ and $V$ are bounded below, 
we find that there exists $\sigma >0$ such that for all $n$,
\[
\int\limits_{\ge_n^{-1}\tilde{\L}} |v_n|^{p+1} \geq \sigma >0.
\]
Hence
\[
t_n^2 \int\limits_{\tilde{\L}} \ge_n^2 |\nabla u_n|^2 + |u_n|^2 \geq \ge_n^N \sigma',
\]
with $\sigma'>0$. Combining this with \eqref{eq:1.24} we see that
\begin{equation}\label{eq:1.23}
t_n \geq \sigma'' >0 \hbox{\quad for all $n\in\mathbb{N}$},
\end{equation}
with $\sigma''>0$. Now, by the definition of $t_n$, we have that
\begin{equation}\label{eq:1.27}
\tilde{E}_n (t_n u_n) \geq \inf_{u\in H^1(\tilde{\L})} \sup_{t>0} \tilde{E}_n (tu)=:b_n.
\end{equation}
But it follows from \cite{dPF3}, Lemma 1.3 with obvious modifications, that
\[
\lim_{n\to\infty} \ge_n^{-N}b_n = c_0,
\]
and this, together with \eqref{eq:1.22}, \eqref{eq:1.23} and \eqref{eq:1.27}, easily 
implies the validity of Claim \eqref{claim1}.

We now proceed to prove Lemma \ref{le:beta}. Set $v_n \colon x \mapsto u_n (\ge_n x)$. We claim that
\begin{equation}\label{claim2}
\sup_{t>0} I_{z_0} (t v_n) \leq c_0 + o(1).
\end{equation}
To see this, we recall that $\{v_n\}$ is bounded in the $H^1$ norm. Since
\[
\int\limits_{\ge_n^{-1}\Omega} \langle J(\ge_n x)\nabla v_n \mid \nabla v_n \rangle + V(\ge_n x) |v_n|^2 = \int\limits_{\ge_n^{-1}\Omega}g(\ge_n x,v_n)v_n \leq \int\limits_{\ge_n^{-1}\Omega} f(v_n)v_n,
\]
similar arguments as those above show that
\[
\int\limits_{\ge_n^{-1}\Omega} |v_n|^{p+1} \geq \sigma >0.
\]
Hence, by the first lemma of Concentration--Compactness (see Lemma I.1 in \cite{Lions}), there is a sequence $B_n$ of balls of radius one such that
\begin{equation}\label{eq:1.33}
\int_{B_n} |v_n|^{2} \geq \sigma >0.
\end{equation}
We now select $t_n >0$ such that $I_{z_0}(t_n v_n)=\sup_{t>0} I_{z_0} (tv_n)$. Since $\{v_n\}$ is bounded in $H^1$ norm we get
\[
C t_n^2 - \int\limits_{\ge_n^{-1}\Omega} F(t_n v_n) \geq I_{z_0} (t_n v_n) \geq c_0.
\]
But from assumption {\bf (f3)} we see that $F(u) \geq C u^\theta$, so that
\[
t_n^{\t -2} \int\limits_{\ge_n^{-1}\Omega} |v_n|^\t \leq C.
\]
This and \eqref{eq:1.33} imply that $\{t_n\}$ is bounded. Therefore from \eqref{claim1} we deduce
\begin{equation}\label{eq:tv}
\lim_{n\to\infty}\int\limits_{\RN \setminus \ge_n^{-1} \tilde{\L}} |t_n v_n|^2 =0.
\end{equation}
From the properties of $z_0$ we easily get (here we use \textbf{(J1)} to get rid of the contribution of $J$ both inside and outside $\tilde{\L}$)
\[
c_0 + o(1) \geq \ge_n^{-N} E^{\ge_n}(t_n u_n) \geq I_{z_0}(t_n v_n) - 
\frac{t_n^2}{2} \int\limits_{\RN \setminus \ge_n^{-1} \tilde{\L}} V(z_0)|v_n|^2.
\]
and so we get \eqref{claim2}.

If we set $w_n = t_n v_n$, we see that $\{w_n\}$ is a minimizing sequence for $I_{z_0}$ constrained to its Nehary manifold $\Ne_{z_0}$. By a straightforward application of the Ekeland variational principle, we can build a Palais--Smale sequence $\{\tilde{w_n}\}$ of $I_{z_0}$ such that $\tilde{w_n}-w_n \to 0$ strongly in $H^1$. Thus there exists a sequence of points $\{z_n\}$ such that $\{w_n (\cdot + z_n)\}$ converges strongly to a positive critical point $w_\infty$ of $I_{z_0}$. Let $\bar y_n=\ge_n z_n$. If 
$\liminf_{n\to\infty} \operatorname{dist}(\bar y_n, \L) >0$ then, since we can take $\tilde{\L}$ 
sufficiently small, we have also $\liminf_{n\to\infty} \operatorname{dist}(\bar y_n, \tilde{\L}) >0$ 
and so from \eqref{eq:tv} we get
\begin{multline*}
o(1)= \int\limits_{\RN \setminus \ge_n^{-1}\tilde{\L}} |t_n v_n|^2 =
\int\limits_{\RN \setminus \ge_n^{-1}\tilde{\L}} |w_n|^2=
\\
= \int\limits_{\RN \setminus \ge_n^{-1}(\tilde{\L}-\bar y_n)} |w_n (\cdot + z_n)|^2 = \int_\RN |w_\infty|^2 + o(1),
\end{multline*}
which contradicts the positivity of $w_\infty$. Hence we may assume that $\bar y_n \to \bar y \in \L$. 
But then
\[
\beta (u_n) = \frac{\int_\RN \chi (\ge_n x +  \bar y_n)|w_n(x+z_n)|^2\, dx}{\int_\RN |w_n (x+z_n)|^2\, dx}
\to \bar y \in \L,
\]
against our (absurd) assumptions $\beta (u_n) \notin \tilde{\L}$.
\end{proof}

\vskip13pt

\begin{proof3}
By Lemma \ref{le:c0}, the map $\xi \mapsto \Phi_\ge(\xi)$ sends $M$ into $\tilde{\Ne}^\ge$. 
Moreover, by Lemma \ref{le:beta} we know that $\beta (\tilde{\Ne}^\ge)\subset \tilde{\L}$. 
Then the map $\xi \mapsto \beta \circ \Phi_\ge(\xi)$ is homotopic to the inclusion 
$j\colon M \to \tilde{\L}$, for any $\ge$ sufficiently small. We now combine Theorem 
\ref{th:cat} with Lemma \ref{le:cat} to get that $E^\ge$ has at least $\cat(M,\tilde{\L})=\cat(M,\L)$ 
critical points on the manifold $\Ne^\ge$. The verification that each one of these 
critical points is actually a solutions of \eqref{eq:E} follows again from Section \ref{sec:pen}, 
once we recall that the main formula in Lemma \ref{le:dis} holds true for each one of the 
critical points just found by definition of $\tilde{\Ne}^\ge$. This completes the proof.
\end{proof3}

\subsection{Using perturbation method}

Let us introduce a topological invariant related to Conley theory.

\begin{definition}
Let $M$ be a subset of $\RN$, $M\ne \emptyset$.
The {\sl cup long} $l(M)$ of $M$ is defined by
$$
l(M)=1+\sup\{k\in \N: \exists\, \a_{1},\ldots,\a_{k}\in
\check{H}^{*}(M)\setminus 1, \,\a_{1}\cup\ldots\cup\a_{k}\ne 0\}.
$$
If no such class exists, we set $l(M)=1$. Here $\check{H}^{*}(M)$ is the
Alexander
cohomology of $M$ with real coefficients and $\cup$ denotes the cup product.
\end{definition}

Let us recall Theorem 6.4 in Chapter II of \cite{C}.

\begin{theorem}\label{th:Chang}
Let $h \in C^2(\RN)$ and let $M\subset \RN$ be a smooth compact nondegenerate manifold of 
critical points of $h$. Let $U$ be a neighborhood of $M$ and let $l \in C^1(\RN)$. 
Then, if $\|h-l\|_{C^1(\bar{U})}$ is sufficiently small, the function $l$ possesses 
at least $l(M)$ critical points in $U$.
\end{theorem}

Let us suppose that $\G$ has a smooth manifold of critical points $M$. We say that $M$
is nondegenerate (for $\G$) if every $x\in M$ is a nondegenerate
critical point of $\G_{|M^{\perp}}$. The Morse index of $M$ is, by definition,
the Morse index of any $x\in M$, as critical point of $\G_{|M^{\perp}}$.

We now can state our multiplicity result.

\begin{theorem}\label{th:main}
Let {\bf (V-V1)} and {\bf (J-J1)} hold and suppose $\G$ has
a nondegenerate smooth manifold of critical points $M$.
Then for $\eps>0$ small, (\ref{eq:E1}) has at least $l(M)$ solutions that
concentrate near points of $M$.
\end{theorem}

\begin{proof}
First of all, we fix $\overline{\xi}$ in such a way that
$|x|<\overline{\xi}$ for all $x\in M$.  We will apply the
finite dimensional procedure with such $\overline{\xi}$ fixed. 

In order to use Theorem \ref{th:Chang}, we set $h(\xi)=C_1 \G(\xi)$ and $l(\xi) = \Phi_\eps(\xi /\eps)$. 
Fix a $\delta$-neighborhood $M_\delta$ of $M$ such that $M_\delta \subset \{ |x|<\overline{\xi}\}$ 
and the only critical points of $\G$ in $M_\delta$ are those in $M$. We will take $U=M_\delta$. 

By \eqref{eq:Phi} and \eqref{eq:DF},  $\Phi_\eps(\cdot /\eps)$ converges to $C_1 \G(\cdot)$ 
in $C^1(\bar{U})$ and so, by Theorem \ref{th:Chang} we have at least $l(M)$ critical points of $l$ 
provided $\eps$ sufficiently small. The concentration statement follows as in \cite{AMS}.
\end{proof}

Clearly, this theorem shows that there is no essential difficulty in dealing with local maxima of $\Gamma$ 
instead of minima. Moreover, when we deal with local minima (resp. maxima) of $\G$, the
preceding results can be improved because the number of positive solutions of \eqref{eq:E1} 
can be estimated by means of the category and $M$ does not need to be a manifold.

\begin{theorem}\label{th:CL}
Let {\bf (V-V1)} and {\bf (J-J1)} hold and suppose $\G$ has
a compact set $X$ where $\G$ achieves a strict local minimum (resp. maximum), 
in the sense that there exists $\delta>0$ and a $\d$-neighborhood $X_\d$ of $X$ such that
\[
b:=\inf\{\G(x):x\in \partial X_{\d}\}>a:= \G_{|X}, \quad
\left({\rm resp. }\; \sup\{\G(x):x\in \partial X_{\d}\}<a\right).
\]

Then  there exists $\eps_{\d}>0$ such that \eqref{eq:E1} has at least $\cat(X,X_\d)$
solutions that concentrate near points of $X_{\d}$, provided $\eps\in
(0,\eps_{\d})$.
\end{theorem}

\begin{proof}
We will treat only the case of minima, being the other one similar.
Fix again $\overline{\xi}$ in such a way that $X_{\d}$ is contained in
$\{x\in\RN : |x|<\overline{\xi}\}$.  We set
$
X^{\eps}=\{\xi:\eps\xi\in X\}$, $X_{\d}^{\eps}=\{\xi:\eps\xi\in X_{\d}\}$ and
$Y^{\eps}=\{\xi\in X_{\d}^{\eps} :\Phi_{\eps}(\xi)\leq C_{1}(a+b)/2\}$.
By \eqref{eq:Phi} it follows that there exists $\eps_{\d}>0$ such that
\begin{equation}\label{eq:X}
X^{\eps}\subset Y^{\eps}\subset X^{\eps}_{\d},
\end{equation}
provided $\eps\in (0,\eps_{\d})$. Moreover, if $\xi\in \partial X_{\d}^{\eps}$ then
$\G(\eps\xi)\geq b$ and hence
\[
\Phi_{\eps}(\xi)\geq C_{1}\G(\eps\xi) + o_{\eps}(1) \geq C_{1}b +
o_{\eps}(1) .
\]
On the other side, if $\xi\in Y^{\eps}$
then $\Phi_{\eps}(\xi)\leq C_{1}(a+b)/2$.  Hence, for $\eps$ small, 
$Y^{\eps}$ cannot meet $\partial X_{\d}^{\eps}$ 
and this readily implies that $Y^{\eps}$ is compact.
Then $\Phi_{\eps}$ possesses at least $\cat(Y^{\eps},X^{\eps}_{\d})$
critical points in $ X_{\d}$.  Using (\ref{eq:X}) and the properties of
the category one gets
\[
\cat(Y^{\eps},Y^{\eps})\geq \cat(X^{\eps},X^{\eps}_{\d})=\cat(X,X_{\d}),
\]
and the result follows.
\end{proof}

\begin{remark}
Let us observe that Theorem \ref{th:ex} is a particular case of Theorem \ref{th:CL}.
\end{remark}

\end{document}